\documentclass[12pt]{article}

\usepackage{amsmath,amssymb,amstext}
\usepackage{amsfonts}

\usepackage[latin1]{inputenc}

\newenvironment{disarray}{\everymath{\displaystyle\everymath{}}\array} {\endarray}
\newtheorem{theorem}{Theorem}
\newtheorem{proposition}{Proposition}

\newtheorem{corollary}{Corollary}
\newtheorem{lemma}{Lemma}

\newtheorem{rk}{Remark}

\newenvironment{proof}{\textbf{Proof}}{\flushright$\blacksquare$\\}
\newcommand\nbOne{{\mathchoice {\rm 1\mskip-4mu l} {\rm 1\mskip-4mu l}
{\rm 1\mskip-4.5mu l} {\rm 1\mskip-5mu l}}}

\title{Multifractal spectra and precise rates of decay in homogeneous fragmentations.}
\author{Nathalie Krell}
\begin{document}
\maketitle

\begin{center}
\it{Laboratoire de Probabilités et Modèles Aléatoires,\\
 Universit\'e Pierre et Marie Curie,\\ 175 rue du Chevaleret,
75013 Paris, France.}
\end{center}

\vspace*{0.8cm}

\begin{abstract}
We consider a mass-conservative fragmentation of the unit
interval. Motivated by a result of Berestycki \cite{be}, the main purpose of this work is to specify the
Hausdorff dimension of the set of locations having exactly an
exponential decay. The study relies on an additive martingale
which arises naturally in this setting, and a class of L\'evy
processes constrained to stay in a finite interval.

\end{abstract}


\bigskip
\noindent{\bf A.M.S. Classification. } 60G09, 60J25, 60J80.

\bigskip
\noindent{\bf e-mail. } nathalie.krell@upmc.fr

\vspace*{0.8cm}

\vspace*{0.8cm}

\section{Introduction.}
Fragmentation appears in a wide range of phenomena in science and
technology, such as degradation of polymers, colloids, droplets,
rocks,... See the proceedings \cite{beycampef} for some applications
in physics, for  example \cite{kramaj} for computer science,
\cite{BertoinMartinez} for mineral crushing,  and works quoted in
\cite{be}  for some further references. This work is a contribution
to the study of the rates of decay of fragments. More precisely, our
aim is to investigate the set of locations which have an exact
exponential decay  (see~ (\ref{Gvab}) below for a precise
definition).

 Roughly a homogeneous fragmentation of intervals $F  (t)$ can
 be seen as a family of nested open sets in $ (0,1)$ such that each interval component
  is spill independently of the others, independently
of the way that spill before, and with the same law as that of the
initial fragmentation  (up to spatial rescaling). We will suppose
that no loss of mass occurs during the process.

Let $x\in (0,1)$ and $I_{x}  (t)$ be the interval component of the
fragmentation $F  (t)$ which contains $x$, and $|I_{x}  (t)|$ its
length.  Bertoin showed in \cite{ber2} that if $V$ is a uniform
random variable on $[0,1]$ which is independent of the
fragmentation, then $ \xi (t):=-\log |I_{V}  (t)|$
  is a subordinator entirely determined by the
  fragmentation characteristics. By the SLLN  for a subordinator, there exists $v_{typ}$ such that $\frac{\xi (t)}{t}\rightarrow
v_{typ} \; a.s.$, which means that $|I_{V}  (t)|\approx
e^{-v_{typ}t}$. Berestycki \cite{be}  computed the Hausdorff
dimension of the set
$$G_{v}:=\left\{x\in  (0,1) \ : \ \lim_{t\rightarrow\infty}
\frac{1}{t}\log |I_{x}  (t)|=-v\right\}$$ for all $v>0$. In this
article we shall rather consider for some $0<a<b$ the set

  \begin{equation} G_{  (v,a,b)}:=\left\{x\in  (0,1)\ : \ a
 \leq \liminf _{t\rightarrow \infty} e^{vt}|I_{x}  (t)|\leq \limsup _{t\rightarrow \infty} e^{vt}|I_{x}  (t)|\leq
b \right\}.
 \label{Gvab}
 \end{equation}

Our goal is to compute  the Hausdorff dimension of the set $G_{
(v,a,b)}$.  Our approach relies on some results on L\'evy processes
constrained to stay in a given interval.

Firstly we will recall background on fragmentations and L\'evy
processes. Secondly we will consider an additive martingale $M$
which is naturally associated to the problem and obtain a criterion
for uniform integrability. This is used in
Section~\ref{limittheorem} to derive some limit theorems which may
be of independent interest  (see Engl\"ander and Kyprianou
\cite{engkyp} for a related approach in the setting of spatial
branching processes). Finally we will compute the Hausdorff
dimension of $G_{  (v,a,b)}$ in Section~\ref{sechausdorffdim}.

\section{Preliminaries.}

\subsection{Definition of fragmentation.}\label{iii}

 We will recall some
 facts about homogeneous interval fragmentations,
which are mostly lifted from \cite{be} , \cite{ber2} and
\cite{ber3}. More precisely, we will consider fragmentations defined
on the space $\mathcal{U}$ of open subsets of $ (0,1)$.  We shall
use the fact that every element $U$ of $\mathcal{U}$ has an interval
decomposition, i.e. there exists a collection of disjoint open
intervals $  (J_{i})_{i\in I}$, where the set of indices $I$ can be
finite or countable, such that $U=\cup_{i\in I}J_{i}$. Each interval
component is viewed as a fragment.

A homogeneous interval fragmentation is a Markov process with
 values in  the space $\mathcal{U}$   which enjoys two keys
 properties. First the branching property: different fragments have independent evolutions. Second, the homogeneity
 property: up to an obvious spacial rescaling, the law of the fragment process does not depend on the initial length of the interval.

 Specifically, if $ \mathbb{P}$ stands for the law of the interval
 fragmentation $F$ started from $F  (0)= (0,1)$, then for $s,t\geq 0$
 conditionally on the open set $F  (t)=\underset{i\in I}{\cup}\ J_{i}  (t)$,
 the interval fragmentation $F  (t+s)$ has the same law as
$F^{1}  (s)\cup F^{2}  (s)\cup ...$
 where for each $i$, $F^{i}  (s)$ is a subset of $J_{i}  (t)$ and has
 the same distribution as the image of $F  (s)$ by the homothetic
 map $ (0,1)\rightarrow J_{i}  (t)$.

\subsection{Poissonian construction of the fragmentation.}\label{pre}
Recall that  $\mathcal{U}$ denotes the space of  open subsets of $
(0,1)$, and set $ \mathbf{1}= (0,1)$. For $U\in \mathcal{U}$,
$$|U|^{\downarrow}:=  (u_{1},u_{2},...)$$ will be the decreasing
sequence of the interval component lengths of $U$. For $U=
(a_{1},b_{1})\in\mathcal{U}$, we define the affine transformation
$g_{U}: (0,1)\rightarrow U$ given by $g_{U}
 (x)=a_{1}+x  (b_{1}-a_{1})$.

In this article we will only consider proper fragmentations  (which
means that the Lebesgue measure of $F  (t)$ is equal to 1). In this
case, Basdevant \cite{ba} has shown that the law of the interval
fragmentation $
 F $ is completely characterized by the so-called dislocation
measure $\nu$  (corresponding to the jump-component of the process)
which is a measure on $\mathcal{U}$ which fulfills the conditions
$$\nu  (\mathbf{1})=0,$$

  \begin{equation}\int_{\mathcal{U}}  (1-u_{1})\nu  (dU)<\infty,
 \label{mesuredelevy}
 \end{equation}
and
$$\sum_{i=1}^{\infty}
u_{i}=1\ \ \ \ \ for\ \nu-almost\ every \ U\ \in \mathcal{U}.$$ This
last assumption is imposed by the hypothesis of length-conservation
and means that when a sudden dislocation occurs, the total length of
the intervals is unchanged. Specialists will notice that the erosion
rates of the fragmentation $c_{r}$ and $c_{l}$ are here equal to 0
for the same reason.

  We now recall the interpretation of
   sudden dislocations of  the fragmentation process in terms
of atoms of a Poisson point process  (see \cite{ba}, \cite{be2}).
 Let $ \nu$ be a dislocation
measure fulfilling the preceding conditions. Let $ K =  (  ( \Delta
 ( t ) , k  ( t ) ), t \geq
 0 )$ be a Poisson point
process with values in $\mathcal{U}\times \mathbb{N}$, and with
intensity measure $\nu\otimes\sharp$, where $\sharp$ is the counting
measure on $\mathbb{N}$. As in \cite{be2}, we can construct a unique
$ \mathcal{U}$-valued process $F=  (F  (t),t\geq 0)$ started from $
(0,1)$, with paths that jump only for times $t\geq 0$ at which a
point $  (\Delta  (t),k  (t))$ occurs, and then $ F  (t)$ is
obtained by replacing the $k  ( t)$-interval $ J_{k
 (t)}  (t-)$ by $ g_{J_{k  (t)}  ( t-)}  (\Delta  (t))$. This point of
view will be used in Section~\ref{sec}.

Some information about the dislocation measure $\nu$ and therefore
about the distribution of the homogeneous fragmentation $F$ is
contained in the function:

 \begin{equation}\kappa  (q)
:=\int_{\mathcal{U}}\left  ( 1-\sum_{j=1}^{\infty}u_{j}^{q+1}\right)
\nu  (dU)\ \ \ \ \forall q>\underline{p} \label{kappa}
\end{equation}
with $\underline{p}$ the smallest real number for which $\kappa$
remains finite :
$$\underline{p}:=\inf\left\{p\in\mathbb{R}:\ \int_{\mathcal{U}}\sum_{j=2}^{\infty}u_{j}^{p+1}
\nu  (dU)<\infty \right\}.$$ We have that $-1\leq \underline{p}\leq
0$  (because $\int_{\mathcal{U}}  (1-u_{1})\nu  (dU)<\infty$ and
$\sum_{i=1}^{\infty} u_{i}=1$ for $ \nu$-almost every $U\ \in
\mathcal{U}$).

This point of view is the same as in \cite{be}  and \cite{ber2},
which deal with ranked fragmentation instead of interval
fragmentation. In the latter the space $ \mathcal{U}$ is replaced by
the space of mass partitions
$$\mathcal{S^{\downarrow}}:=\left\{x=  (x_{1},x_{2},...)\ |\
 x_{1}\geq x_{2}\geq ...\geq 0\ , \sum_{i=1}^{\infty}
x_{i}\leq 1 \right\}.$$For the precise link between these two
fragmentations see \cite{ba}.

\subsection{An important
subordinator.}\label{subordinateur}

Let $x\in (0,1)$ and $I_{x}  (t)$ be the interval component of the
random open set $F  (t)$ which contains $x$, and $|I_{x}  (t)|$ its
length. Let $V$ be a uniform random variable on $[0,1]$ which is
independent of the fragmentation.

Bertoin showed in \cite{ber2} that
\begin{equation} \xi (t):=-\log |I_{V}  (t)|,\ \ \ t\geq 0 \label{xi},
 \end{equation}
  is a subordinator, with
Laplace exponent $\kappa  (q)$ defined in~ (\ref{kappa})  (i.e.
$\mathbb{E}  (e^{-\lambda\xi (t)})=e^{-t\kappa  (\lambda)}$ for all
$\lambda >\underline{p}$). In order to interpret this as a
L\'evy-Khintchine formula, we introduce the measure
$$L  (dx):=e^{-x}\sum_{j=1}^{\infty}\nu  (-\log u_{j}\in dx),\ \ \ \ \
\ x\in  (0,\infty ).$$ It is easy to check that $\int \min  (1, x) L
 (dx)<\infty$, thus $L$ is the L\'evy measure of a subordinator,
and we can check that $\kappa  (q)=\int_{  (0,\infty )}\left
(1-e^{qx}\right) L  (dx).$

In this article we shall consider the L\'evy process $Y_{t}=vt-\xi
(t)$. In order to apply certain results to this process, we will
need to assume that its one-dimensional distributions are absolutely
continuous. Let $L^{ac}$ be the absolutely continuous part of the
measure $L$.   Tucker  has shown in \cite{tuc} that

\begin{equation} \int_{\mathbb{R}_{+}}\frac{1}{1+x^{2}} L^{ac}  (dx)=\infty,\label{lac}
 \end{equation}
ensures the absolute continuity of  one-dimensional distribution of
the L\'evy
 process evaluated at any $t>0$.
As  $\int \min  (1, x) L
 (dx)<\infty$, the condition  (\ref{lac}) is equivalent to :
 \begin{equation}
L^{ac}  ([0,\epsilon))=\infty \ \ \ \text{for any }\ \ \ \epsilon>0.
 \label{lac2}
 \end{equation}
 Let $\nu_{1}$ be the image of the measure $\nu$ by the map $U\rightarrow u_{1}$  (recall that $u_{1}$ is the length of the longest interval component of the open set $U$) and $\nu_{1}^{ac}$ be the
absolutely continuous part of the measure $\nu_{1}$. Throughout this
work we will make the following assumption, which is easily seen to
imply  (\ref{lac2})  (in fact we can even show that the two are
equivalent):
 \begin{equation}
\nu_{1}^{ac}  ([0,\epsilon))=\infty \ \ \ \text{for any }\ \ \
\epsilon>0.
 \label{nuinf}
 \end{equation}
In the next subsection, we will give some results about L\'evy
processes that will be needed in the sequel, and apply for
$Y_{t}=vt-\xi (t)$.

\subsection{An estimate for completely asymmetric L\'evy processes.}

For the next sections, we will need some technical notions about
completely asymmetric L\'evy processes. Therefore we recall some
facts mostly lifted from \cite{BertoinLevy} and \cite{ber1}. Let $Y=
( Y_{t})_{t\geq0}$ be a L\'evy process with no positive jumps and $
  (\mathcal{E}_{t})_{t\geq 0}$ the natural filtration associated to
$  (Y_{t})_{t\geq0}$. The case where $Y$ is the negative of a
subordinator
 is degenerate for our purpose and
therefore will be implicitly excluded in the rest of the article.
The law of the L\'evy process started at $x\in \mathbb{R}$ will be
denoted by $\mathbf{P}_{x}$  (so bold symbols $\mathbf{P}$ and
$\mathbf{E}$ refer to the L\'evy process while $\mathbb{P}$ and
$\mathbb{E}$ refer to the fragmentation), its Laplace transform is
given by
$$\mathbf{E}_{0}  (e^{\lambda Y_{t}})=e^{t\psi  (\lambda)}, \ \
\lambda,\ t \geq0,$$ where $\psi:
\mathbb{R}_{+}\rightarrow\mathbb{R}$ is called the Laplace exponent.

Let $\phi:\mathbb{R}_{+}\rightarrow\mathbb{R}_{+}$ be the right
inverse of $\psi$  (which exists because
$\psi:\mathbb{R}_{+}\rightarrow\mathbb{R}$ is convex with
$\lim_{t\rightarrow\infty}\psi  (\lambda)=\infty$), i.e. $\psi  (
\phi   (\lambda))=\lambda\ \ \ \forall \lambda\geq0$.

Let us recall some important features on the two-sided exit problem
(which is completely solved in \cite{ber1}). For $\beta>0$ we denote
the first exit time from $ (0,\beta )$ by

\begin{equation}
T_{\beta}=\inf\{t : \ Y_{t} \notin  (0,\beta )\}. \label{ta}
 \end{equation}
Let $W:\mathbb{R}_{+}\rightarrow\mathbb{R}_{+}$ be the scale
function, that is  the unique
 continuous function  with Laplace transform:
$$\int_{0}^{\infty}e^{-\lambda x} W  (x) dx=\frac{1}{\psi  (\lambda)}\
\ \ , \ \lambda>\phi  (0).$$

For $q\in\mathbb{R}$, let $W^{   (
q)}:\mathbb{R}_{+}\rightarrow\mathbb{R}_{+}$ be the
 continuous function such that for every $x\in\mathbb{R}_{+}$
 $$W^{  (q)}  (x):=\sum_{k=0}^{\infty}q^{k}W^{*k+1}  (x),$$
 where $W^{*n}=W*...*W$ denotes the $n$th convolution power of
 the function W  (for more details about this see \cite{BertoinLevy} or
 \cite{ber1}). So that
$$\int_{0}^{\infty}e^{-\lambda x} W^{ (q)}  (x) dx=\frac{1}{\psi  (\lambda)-q}\
\ \ , \ \lambda>\phi  (q).$$

The next statement is about the
 asymptotic behavior of the  L\'evy process killed when it exits $ (0,\beta )$  (point 1 and 2), which is taken from \cite{ber1}, and about   the L\'evy process conditioned to remain in $ (0,\beta )$  (point 3, 4 and 5), which is taken from Theorem 3.1
 ( ii)
and Proposition 5.1  (i) and  (ii) in \cite{lam} :

\begin{theorem}\label{theo11}
 Let us define the transition probabilities
 $$ P_{t}  (x,A):=\mathbf{P}_{x}  (Y_{t}\in A , t<T_{\beta})\
 \hbox{for} \ x\in (0,\beta )\ \hbox{and} \ A\in\mathcal{B}  ( (0,\beta )),$$ and the
 critical value
 \begin{equation}
  \rho_{\beta}:=\inf\{ q\geq 0\ ;\
 W^{  (-q)}  (\beta)=0\}\label{rho},
 \end{equation}
 Suppose that the one-dimensional distributions of the L\'evy
 process are absolutely continuous. Then the following holds:
 \begin{enumerate}
\item $\rho_{\beta}\in  (0,\infty )$
 and the function $W^{  (-\rho_{\beta})}$ is strictly positive on $ (0,\beta )$

 \item Let $\Pi  (dx):=W^{  (-\rho_{\beta})}  (\beta-x)dx$.
For every $x\in (0,\beta )$:
$$\lim_{t\rightarrow\infty} e^{\rho_{\beta} t} P_{t}  (x,.) = c
W^{  (-\rho_{\beta})}  (x) \Pi  (.)$$ in the sense of weak
convergence, where
$$c:=\left  (\int_{0}^{\beta}W^{  (-\rho_{\beta})}  (y) W^{
  (-\rho_{\beta})}  (\beta-y)dy\right)^{-1}.$$
 \item The
process
\begin{equation}
D_{t}:=e^{\rho_{\beta} t}\ \mathbf{1}_{\{t<T_{\beta}\}}\ \frac{W^{
  (-\rho_{\beta})}  (Y_{t})}{ W^{  (-\rho_{\beta})}  (x)}\label{dt}
\end{equation}
 is a $  (\mathbf{P}_{x},  (\mathcal{E}_{t}))$-martingale.

\item The mapping $  (x,q)\mapsto W^{  (q)}  (x)$ is of class $\mathcal{C}^{1}$ on $ (0,\infty )\times  (-\infty,\infty).$
 \item The mapping
$\beta\mapsto\rho_{\beta}=\inf\{q>0 : W^{  (-q)}  (\beta)=0\}$ is
strictly decreasing and of class $ \mathcal{C}^{1}$ on $ (0,\infty
)$.

 \end{enumerate}
\end{theorem}

\begin{rk}
 The definition of $\rho_{\beta}$ is of course complicated, however in the  simple case  when $Y$ is a standard Brownian motion, we have:
$$\rho_{\beta}=\pi^{2}/\beta^{2}\ \ \  \hbox{and }\ \ \  W^{ (-\rho_{\beta})} (x)=\frac{\beta}{\pi}\sin\left (\frac{\pi}{\beta}x\right).$$
In the case where $Y$ is a standard stable process, the mapping of
$\beta \rightarrow \rho_{\beta}$ is depicted in \cite{ber1996}. We
also point at the more explicit lower bound  (see Lemma 5 in
\cite{ber1}):
$$\rho_{a}\geq 1/W (a),$$
Another lower bound will be given in Remark \ref{rk sur c  (v)}
below.
\end{rk}

\begin{rk}
 The formula for the constant $c$ in part 2. of Theorem \ref{theo11} stems from the relation
$$ e^{\rho_{\beta} t} \frac{W^{  (-\rho_{\beta})}  (y)}
{W^{  (-\rho_{\beta})}  (x)} P_{t}  (x,dy)
\underset{t\rightarrow\infty}{\sim} c  W^{  (-\rho_{\beta})}   (
\beta-y)W^{  (-\rho_{\beta})}  (y)dy.$$ Integrating over $ (0,\beta
)$ and using the fact that $D_{t}$  is a martingale yields the given
expression.

\end{rk}

We also refer to the recent article of T. Chan and A. Kyprianou
\cite{ChaKyp} for further properties of $W^{  (-\rho_{\beta})}$.

Now we have recalled the background that  is needed to solve our
problem.

\section{An additive martingale. }\label{sec}

 Now we turn our attention to the main purpose of this article and consider a
homogeneous interval fragmentation $  (F  (t),t\geq 0)$ and some
real numbers  $v>0$ and $0<a<b$. We are interested in the asymptotic
set:
$$G_{  (v,a,b)}=\left\{x\in  (0,1)\ : \ a
 \leq \liminf _{t\rightarrow \infty} e^{vt}|I_{x}  (t)|\leq \limsup _{t\rightarrow \infty} e^{vt}|I_{x}  (t)|\leq
b \right\},
$$ with $|I_{x}  (t)|$ the length
of the interval component of $F  (t)$ which contains $x$.

In order to do that, we will have to consider first the non
asymptotic set:
$$\Lambda_{  (v,a,b)}=\left\{x\in  (0,1)\ : \  ae^{-vt}< |I_{x}  (t)|< b
e^{-vt} \ \forall t\geq 0\right\},$$ for  $0<a<1<b$.

\textbf{In this section and in the next we will assume that
$0<a<1<b$.}

We introduce some notation, that we will need in the rest of the
article: define the set of the ``good'' intervals at time $t$ as
\begin{equation}
G  (t):=\{I_{x}  (t): \ x\in (0,1) \ \ \hbox{and} \  \ a e^{-vs}<
|I_{x}  ( s)|< b e^{-vs} \ \ \forall\ s\leq t\}.\label{good}
\end{equation}

Let $  (\mathcal{F}_{t})_{t\geq 0}$ be the natural filtration of the
interval fragmentation $  (F  (t),t\geq 0)$. Let $
(\mathcal{G}_{t})_{t\geq 0}$ be the enlarged filtration defined by $
\mathcal{G}_{t}= \mathcal{F}_{t}\vee\sigma (I_{V} (t))$ where $V$ is
a uniform variable independent of the fragmentation). We can remark
that for all $t$ we have $
\mathcal{G}_{t}\subsetneq\mathcal{F}_{t}\vee \sigma\{V\}$, and $
\mathcal{G}_{\infty}=\mathcal{F}_{\infty}\vee \sigma\{V\}$.

We recall that $\xi (t)=-\log |I_{V}  (t)|$ is a subordinator. More
precisely we are interested in the  L\'evy process with no positive
jump $Y_{t}:=vt-\xi (t)+\log  (1/a)$, and use the results of
preceding subsection for this L\'evy process. We remark that its
Laplace exponent $\psi  (\lambda)$ is equal to $v\lambda-\kappa  (
\lambda)$, with $\kappa$ defined in Subsection~\ref{subordinateur}.
 Since we have supposed~ (\ref{nuinf}), the one-dimensional distributions of the L\'evy
 process $Y_{t}$ are absolutely continuous and we
can apply Theorem~\ref{theo11}.

 For this L\'evy process
$Y$  let $$T:=T_{\log  (b/a)}$$  and
$$\rho:=\rho_{\log  (b/ a)},$$ where $T_{\beta}$ is defined in~ (\ref{ta}) and $\rho_{\beta}$ is defined in~ (\ref{rho}). We stress
 that $\rho$ depends on $v$, $a$, $b$ and $\kappa$.

 To simplify the notation, let also
 $$h  (t):= W^{  (-\rho)}  (t-\log
a)\mathbf{1}_{\{t\in  (\log a ,\log b )\}}$$ for all
$t\in\mathbb{R},$ and $h (-\infty)=0$.

By rewriting~ (\ref{dt}) with the new notation  we get a $  (
\mathcal{G}_{t})$-martingale
$$D_{t}=e^{\rho t}\ \mathbf{1}_{\{t<T\}}\ \frac{h
  (vt+\log |I_{V}  (t)|)}{h  (0)},\ \ \  \ \  t\geq 0.$$
If $I$ is an interval component of $F  (t)$, we define the
``killed'' interval $I^{\dag}$ by $I^{\dag}=I$ if $I$ is good  (i.e.
$I\in G
 (t)$ with $ G  ( t)$ defined in~ (\ref{good})), else by
$I^{\dag}=\emptyset$. Projecting the martingale $D_{t}$ on the
sub-filtration $   ( \mathcal{F}_{t})_{t\geq 0}$, we obtain an
additive martingale
$$
M_{t}:=\frac{e^{\rho t}}{h  (0)}\int_{0}^{1} h  (vt+\log
|I_{x}^{\dag}  (t)|)\ dx\ ,\ \ \  \ \  t\geq 0.$$ We notice that if
$y\in I_{x}  (t)$, then $I_{y}  (t)=I_{x}  (t)$. Now we will
consider the interval decomposition $  (J_{i}  (t), J_{2}  (t),...)$
of the open $F  (t)$  (see subsection \ref{iii}). We can rewrite
$M_{t}$ as:
\begin{eqnarray} M_{t}=\frac{e^{\rho t}}{h  (0)} \sum_{i\in\mathbb{N}}
h\left  (vt+\log |J_{i}^{\dag}  ( t)|\right)|J_{i}^{\dag}
 (t)|.\label{minf}
\end{eqnarray}
We will use this expression in the rest of the article.

Finally, let  the  absorption time of $M_{t}$ at $0$ be
 $$\begin{disarray}{rcl}
\zeta&:=&\inf\{t :M_{t}=0\}\\&&\\&=&\inf\{t: G
(t)=\emptyset\},\end{disarray}$$ with the convention
$\inf{\emptyset}=\infty$.

Our first result is:

\begin{theorem}\label{theo2}
In the previous notation, with the assumptions~
 (\ref{nuinf}) and  if $v>\rho$ holds, then:
 \begin{enumerate}
\item The martingale $M_{t}$ is bounded in $\mathrm{L}^{2}  (\mathbb{P})$.
 \item Conditionally on $\zeta=\infty$, we have: $ \lim_{t\rightarrow\infty} M_{t}>0 .$
 \end{enumerate}
\end{theorem}

\begin{rk}\label{remarquebplusgrandquedeuxa} We stress that as  $\rho$ depends on $v$, $a$, $b$ and $\kappa$, the condition $v>\rho$ involves implicitly the parameters $a$ and $b$.
In particular it forces $b>2a$, otherwise there would never be more
than one ``good'' interval  (as a fragment of size $x$  will split
into at least  two  different fragments and the smallest one will
have a size at most equal to $x/2$), and as a consequence we would
have $  M_{\infty}=0$ a.s., in contradiction with the  uniform
integrability of $M_{.}$.

\end{rk}

The proof of Theorem \ref{theo2}.1. is given in the appendix.

\bigskip

In order to prove  Theorem \ref{theo2}.2 we will first introduce
some notation, then prove two lemmas, and after we will conclude.

Let $I$ be an interval of $ (0,1)$. The law of the homogeneous
interval fragmentation started at $I$ will be denoted by
$\mathbb{P}_{I}$. We remark that $\mathbb{P}_{I}   (
M_{\infty}=0|\zeta=\infty)$ only depends on the length of $I$.
Therefore we define $$g  (x):=\mathbb{P}_{I}   (
M_{\infty}=0|\zeta=\infty),$$ where $I$ is an interval such that
$|I|=x$. Let $N$ be the integer part of $  (2b-a)/a$. As we assume
$v>\rho$, we have necessarily $b>2a$  (see
Remark~\ref{remarquebplusgrandquedeuxa}), thus $N\geq 2$. Let
$\eta:=
  (b-a)N^{-1}$. We remark that $\eta<a$ and $b-a=N \eta$.
Denote the first time when there are at least two good intervals by
$$T^{F}:=\inf\{t: \sharp G  (t)\geq 2\},$$ with the convention $\inf
\emptyset=\infty$. We notice that $T^{F}$ is an $   (
\mathcal{F}_{t})$ stopping time as $\sharp G  (t)$ is $
\mathcal{F}_{t}$-adapted.

\begin{lemma}\label{lemme1} In the previous notation, supposing that~ (\ref{nuinf}) and
$v>\rho$ hold, we get: for every open interval $I$
$$ \mathbb{P}_{I}  (T^{F}=\infty|\zeta=\infty)=0.$$
\end{lemma}

\begin{proof}
We notice that, as  the martingale $M_{t}$ is not identically 0 and
is uniformly integrable, we have $\mathbb{P}_{I}   (
T^{F}=\infty|\zeta=\infty)<1$  (because $M_{\infty}=0$ when
$T^{F}=\infty$).

Let $I$ be an open interval such that $|I|\in  (a,b)$, $t_{0}:=\log
 ( 2b/a)/v$ and  $\epsilon:=a^{2}/  (2b^{2})$. Thus
$$\begin{disarray}{rcl}|I|  (1-\epsilon)>a/2\geq b e^{-vt_{0}}&
\hbox{and} & |I|\epsilon <b\epsilon\leq ae^{-vt_{0}}\end{disarray}$$
therefore, if the dislocation of $I$ produces at time $t_{0}$ an
interval of length at least $|I| (1-\epsilon)$ then this interval is
too large to be good and the remaining ones are too small to be good
either. As a consequence we have
$$ \mathbb{P}_{I}  (M_{ t_{0}}=0)\geq
\mathbf{P}_{\log |I|}  (e^{-\xi (t_{0})}>e^{-\log |I|}  (
1-\epsilon) )=\mathbf{P}  (\xi (t_{0})<-\log  (1-\epsilon)),$$ by
the homogeneous property of the fragmentation. Moreover since $\xi
(t)$ is a subordinator, we get $p:=\mathbf{P}  ( \xi (t_{0})<-\log
(1-\epsilon))>0$, therefore

\begin{eqnarray}
\mathbb{P}_{I}  (M_{t_{0}}=0)&\geq& p>0.\label{defdep}
\end{eqnarray}
Additionally for every open interval $I$ such that $|I|\in  (a,b)$:
$$ \mathbb{P}_{I}  (\sharp G  (t)=1\ \forall t\leq
t_{0})\leq 1-\mathbb{P}_{I}  (M_{t_{0}}=0)\leq 1-p.$$ Using the
strong Markov property of the fragmentation and~ (\ref{defdep}) we
find by induction that for all $k\in\mathbb{N}$:
$$ \mathbb{P}_{I}  (\sharp G  (t)=1\ \forall t\leq
kt_{0})\leq  (1-p)^{k}.$$ Therefore
$$ \lim_{t\rightarrow\infty}\mathbb{P}_{I}  (\sharp G  (s)=1\ \forall s\leq
t)=0$$ and as a consequence
$$\mathbb{P}_{I}  (T^{F}=\infty|\zeta=\infty)=0.$$
\end{proof}

\begin{lemma}\label{lemme2} In the previous notation, supposing that~ (\ref{nuinf}) and
$v>\rho$ hold, we get:
$$\underset{ a<x<b}{\sup} g  (x)=
\underset{1\leq k\leq N}{\max}g  (a+k\eta),$$where $N=\lfloor
(2b-a)/a \rfloor$ and $\eta= (b-a)/N$.
\end{lemma}

\begin{proof}
  We will prove this lemma by  induction.

  The hypothesis of induction is for $n\leq N$:
  $$\begin{disarray}{lc}
  \left  (H \right)_{n}\ :\ \ \ \ \ \ \ \ \ & \ \ \ \ \sup_{x\in  (a,a+n\eta )} g  (x)=\max_{1\leq
k\leq n} g  (a+\eta k).
 \end{disarray}$$

\medskip

 $\ast$ The case $n=1$: let $I$ be an open interval such that
 $|I|\in (a,a+\eta )$. We work under $ \mathbb{P}_{I}$ conditionally
 on  ``non-extinction"  (which means conditionally on the event
 $\zeta=\infty$). Let $$T^{1}:=\inf\{t\geq 0| \ \exists J  (t)\in
 G  (t)\ : \ e^{vt}|J  (t)|\notin  (a,a+\eta )\},$$
with $G  (t)$ defined in~ (\ref{good}). The random time $T^{1}$ is
an $  (\mathcal{F}_{t})$ stopping times. As the quantity $vt-(-\log
|J
  (t)|)$ creeps upwards with probability equals to 1 and as $J  (t)\in G  (t)$ implies that
$e^{vt}|J  (t)|>a$, we get
$$T^{1}=\inf\{t\geq 0| \ \exists J  (t)\in
 G  (t)\ : \ e^{vt}|J  (t)|=a+\eta\}.$$
Moreover by the choice of $\eta$ we have $a+\eta<2a$, which implies
that there is at most one good interval whose length is always in $
(a,a+\eta )$. Recall from Lemma~\ref{lemme1} that  $ \mathbb{P}_{I}
(T^{F}<\infty| \zeta=\infty)=1$, thus
$$\mathbb{P}_{I}  (T^{1}<\infty| \zeta=\infty)=1.$$
Using the strong Markov property at the stopping times $T^{1}$, we
get $$g  (x)\leq g  (a+\eta)\ \ \ \ ,\ \ x\in (a,a+\eta ),$$ thus $
(H)_{1}$ holds.

\medskip

 $\ast$ The case $n+1$  (with $n+1\leq N$): we suppose that the hypothesis of induction
 holds for all $k \leq n$.

 Let $I$ be an open interval such that
 $|I|\in (a+n\eta,a+  (n+1)\eta)$. We work under $ \mathbb{P}_{I}$ conditionally
 on ``non-extinction". Let $$T^{n}:=\inf\{t\geq 0| \ \exists J  (t)\in
 G  (t)\ : \ e^{vt}|J  (t)|\notin  (a+n\eta,a+  (n+1)\eta )\},$$
with $G  (t)$ defined in~ (\ref{good}). The random time $T^{n}$ is
an $  (\mathcal{F}_{t})$ stopping times. As the quantity $e^{vt}|J
  (t)|$ grows only continuously, we get
$$T^{n}=\inf\{t\geq 0| \ \exists J  (t)\in
 G  (t)\ : \ e^{vt}|J  (t)|=a+  (n+1)\eta\ \hbox{or} \ e^{vt}|J  (t)|\in  (a,a+n\eta]\}.$$
Moreover by the choice of $\eta$ we have $a+\eta<2a$, which implies
that there is at most one good interval which length is always in
$(a+n\eta,a+  (n+1)\eta )$. Additionally by Lemma~\ref{lemme1}, we
get $ \mathbb{P}_{I}  (T^{F}<\infty| \zeta=\infty)=1$, thus
$$\mathbb{P}_{I}  (T^{n}<\infty| \zeta=\infty)=1.$$
Using the strong Markov property at the stopping times $T^{n}$, we
get

 $$\begin{disarray}{lcl}
g  (|I|)&\leq& \max\left  (g  (a+  (n+1)\eta),\sup_{y\in
(a,a+n\eta]}
 g  (y)\right).\end{disarray}$$As this holds for every open interval
 $I$ such that $|I|\in  (a+n\eta,a+  (n+1)\eta )$, by the
hypothesis of induction, we have established
 $ (H)_{n+1}$.
\end{proof}

\begin{proof}[Proof of Theorem \ref{theo2}.2.] With Lemma~\ref{lemme2}, we get that
there exists a integer $k_{0}$ in $[1,N]$ such that $g  (a+\eta
k_{0})=\sup_{x\in  (a,b)} g  (x)$  (if two or more values of $k$,
are possible, we choose the smallest one). Let $x_{0}$ be $a+\eta
k_{0}$.

Additionally, with Lemma~\ref{lemme1}, we get $\mathbb{P}_{
(0,x_{0})}  (T^{F}< \infty|\zeta=\infty)=1$. Using the strong
property of Markov for the stopping times $T^{F}$, and with $n\geq
2$ the random number of good intervals of the fragmentation at time
$T^{F}$ and with $\alpha_{1},..., \alpha_{n}$ the length of those
intervals, we get:
$$g  (x_{0})\leq\mathbb{E}  (g  (\alpha_{1})...g  (\alpha_{n}))\leq
\mathbb{E}  (g  (x_{0})^{n})\leq g  (x_{0})^{2}.$$ As $g  (x_{0})<1$
by the uniformly integrability of $M_{t}$, we get that $g  (
x_{0})=0$ and finally that $g\equiv 0$.
\end{proof}

\section{Limit theorems.}\label{limittheorem}
In this section, we establish two corollaries of
Theorem~\ref{theo2}, which will be useful in the sequel.

Bertoin and Rouault  (Corollary 2 in \cite{berrou}) proved that
\begin{equation}\lim_{t\rightarrow\infty}\frac{1}{t}\log \sharp\{I_{x}  (t):\
ae^{-vt}<|I_{x}  (t)|<be^{-tv} \}=C
(v),\label{eqberrou}\end{equation} where $C  (v):=   (
\Upsilon_{v}+1)v-\kappa  (\Upsilon_{v})$ and $\Upsilon_{v}$ is the
reciprocal of $v$ by $\kappa^{'}$ i.e, $\kappa^{'}   (
\Upsilon_{v})=v$ for $v\in (v_{min},v_{max})$. \footnote{Where
$v_{min}$ is the maximum of the function $p\mapsto \kappa  (p-1)/p$
on $ (\underline{p}+1,\infty )$  and $v_{max}:=\kappa^{'}
 (\underline{p}^{+}) $ (see \cite{be} ).}

Here we deal with the more stringent requirement: $\forall s\leq t,\
|I_{x}  (s)|\in  (ae^{-sv},be^{-sv})$, and the next proposition
gives the rates that we find in that case.

\begin{proposition}\label{theo4}
In the notation of the previous sections, with the assumptions~
 (\ref{nuinf}) and if $v>\rho $ we get that conditionally on
$\zeta=\infty$  (i.e. $M$ is not absorbed at $0$, or in a equivalent
way $\Lambda_{ (v,a,b)}\neq \emptyset$):
\begin{equation}
\lim_{t\rightarrow\infty}\frac{1}{t}\log \sharp G  (t)=\ v-\rho\ \ \
\ \ a.s. \label{eq3}
\end{equation}
\end{proposition}

Before proving this corollary we make the following remark
\begin{rk}\label{rk sur c  (v)}
 It is interesting to compare the  estimate found by Bertoin and
Rouault  and
 the present one
  (of course we have not considered the same set,
nevertheless the two estimates are related). For this we show that
for all $v\in  (v_{min},v_{max})$ and $a$ and $b$ such that
$\rho\geq v_{min}$ we have $C  (v)\geq v-\rho$. In this direction we
use results from \cite{be}  Section 1. Let $\Psi  (p):=p\kappa^{'}
 (p)-\kappa  (p)$ for all $p>0$ with $\kappa^{'}$ the derivative of
$\kappa$  (this function is well defined because of the definition
of $\underline{p}$ in Section 2 and because $\underline{p}\leq 0$).
For every $p>0$, $\Psi^{'}  (p)=p\kappa^{''}  (p)\leq 0$ since
$\kappa$ is concave. As a consequence $\Psi$ is decreasing. With the
definition of $\Upsilon_{v}$, we get that the function $v\in
(v_{min},v_{max})\mapsto \Upsilon_{v}\in\mathbb{R}$ is decreasing,
additionally $\Upsilon_{v_{min}}>0$, therefore the function $v\in
(v_{min},v_{max})\mapsto g
  (\Upsilon_{v})\in \mathbb{R}$ is increasing. Moreover $\Psi   (
\Upsilon_{v})=C  (v)-v$, hence for all $v\in  (v_{min},v_{max})$:
$$C  (v)-v\geq C  (v_{min})-v_{min}=-v_{min}.$$
Additionally as $\rho\geq v_{min}$,  we  finally obtain:
$$  \forall v\in  (v_{min},v_{max})\ \ \ \ C  (v)\geq v-\rho.$$
As a consequence, we have checked that the rate of growth of $\sharp
G  (t)$  (defined in~ (\ref{good})) is lower that of $\sharp \{I_{x}
(t):\ |I_{x}  (t)|\in  (ae^{-tv},be^{-vt})\}$, which was of course
expected.
\end{rk}

\begin{proof}
 In this proof we work conditionally on $\zeta=\infty$  ( i.e $M$
is not absorbed at 0). Applying Theorem~\ref{theo2}, we get $
M_{\infty}>0. $   In order to show that~ (\ref{eq3}) holds, we will
first look at the lower bound of the inequality, and then at the
upper bound.

\medskip

$\bullet$   With the definition  of $M_{t}$ in~ (\ref{minf}), of
 $G  (t)$ and of $J_{i}^{\dag}  (t)$ at the beginning of
Section~\ref{sec} and by the conditioning, there exists $t^{'}>0$
such that for all $t\geq t^{'}$:
$$\begin{disarray}{rcl}
\frac{M_{\infty}}{2}&\leq & \frac{ e^{\rho t}}{h   (
0)}\sum_{i\in\mathbb{N}} h  (vt+\log  (|J_{i}^{\dag}  (t)|))\
|J_{i}^{\dag}  (t)|\leq  \frac{ e^{\rho t}}{h
 (0)}\sum_{i\in\mathbb{N}} C_{4}\ be^{-vt}\ \mathbf{1}_{\{ J_{i}
 (t)\in G  (t)\}},
\end{disarray}$$with $C_{4}$ as maximum of $h  (.)$ on
$[\log a,\log b]$. Hence for all $t\geq t^{'}$ :

$$  \sharp G  (t)\geq\ e^{  (v-\rho)
t} \frac{h  (0)}{2C_{4}b}M_{\infty},$$ and as a consequence,
conditionally on $\zeta=\infty$,
\begin{equation}  \liminf_{t\rightarrow\infty}\frac{1}{t}\log \sharp G  (t)\geq\ v-\rho\label{eq1}.
\end{equation}

$\bullet$ Secondly we will show the converse inequality.

 Let $0< a^{'}<a<1<b<b^{'} $, and
$\rho^{'}:=\rho_{\log  (b^{'}/a^{'})}$.  Denote the set of  ``good''
intervals associated to $a^{'}$ and $b^{'}$ by:
$$G^{'}  (t):=\{I_{x}  (t) : \ \  x\in (0,1) \ \  \hbox{and}\ \  \ |I_{x}  (s)|\ \in \
 (a^{'} e^{-vs},b^{'} e^{-vs} )\ \ \ \forall\ s\leq t\}.$$ Let
$M_{t}^{'}$ be the martingale defined at the beginning of
Section~\ref{sec}  (and denoted there by $M$) associated to $a^{'},
b^{'}$ instead of $a,b$. Plainly, if $M_{t}$ is not absorbed at 0,
then a fortiori $M_{t}^{'}$ is not absorbed at 0 either.
Additionally, since $\log  (b^{'}/a^{'})>\log  (b/a)$, and
$\rho_{.}$ is strictly decreasing   (see Theorem~\ref{theo11}.5), we
get $v>\rho>\rho^{'}$ and we may  apply Theorem~\ref{theo2} for
$a^{'}, b^{'}$ instead of $a,b$. We get $\lim_{t\rightarrow\infty}
M_{t}^{'}=M_{\infty}^{'}>0$.

With the definition~ (\ref{minf}) of $M_{t}$ and with an analogue of
the function $h  (t)$, namely $t\in\mathbb{R}$
$$\varphi  (t):=
W^{  (-\rho^{'})}  (t+\log  (1/a^{'}))\mathbf{1}_{\{t\in  (\log
a^{'},\log b^{'})\}},$$
 we get:
$$\begin{disarray}{rcl}
M_{\infty}^{'}&=&\lim_{t\rightarrow\infty} \frac{ e^{\rho^{'}
t}}{\varphi  (0)}\sum_{i\in\mathbb{N}} \varphi  (vt+\log |J_{i}  (
t)|) |J_{i}  (t)|\ \mathbf{1}_{\{J_{i}  (t)\in G^{'}  (t)\}}.
\end{disarray}$$
Therefore there exists $t^{'}>0$ such that for every $t\geq t^{'}$
$$\begin{disarray}{rcl}
2M_{\infty}^{'}&\geq& \frac{ e^{\rho^{'} t}}{\varphi   (
0)}\sum_{i\in\mathbb{N}} \varphi  (vt+\log |J_{i}  (t)|) |J_{i}
  (t)|\ \mathbf{1}_{\{ J_{i}  (t)\in G^{'}  (t)\}}
\\&&\\&\geq &
 \frac{ e^{\rho^{'}
t}}{\varphi  (0)}\sum_{i\in\mathbb{N}} \varphi  (vt+\log |J_{i}
  (t)|)\
 a^{'} e^{-vt} \ \mathbf{1}_{\{J_{i}  (t)\in G  (t)\}}.
\end{disarray}$$Since $ (ae^{-vt}, be^{-vt})\subsetneq  (a^{'}e^{-vt},
b^{'}e^{-vt})$, we get by Theorem~\ref{theo11}.1, that for all $x\
\in [ \log a, \log b ] $: $\varphi  (x)>0.$ Because $[\log a, \log
b]$ is compact and $\varphi  (.)$ is a continuous function,
$$\underset{x\in[\log a, \log b]}{\inf}\varphi  (x)>0 .$$
Combining this  with $$C_{5}:=2M_{\infty}^{'}\varphi  (0)/\left   (
a^{'}\underset{x\in[\log a, \log b]}{\inf}\varphi   (
x)\right)<\infty,$$ we get for all $t\geq t^{'}$ :
$$C_{5}\geq e^{  (\rho^{'}-v)t} \sum_{i\in\mathbb{N}}
\mathbf{1}_{\{ J_{i}  (t)\in G  (t)\}}$$ and thus
$$C_{5}e^{  (v-\rho^{'})t}\geq \sharp G  (t).$$
Hence for all $a^{'}, b^{'}$ such that $0< a^{'}<a<1<b<b^{'}$:
 $$\limsup_{t\rightarrow\infty}\frac{1}{t}\log\sharp G  (t)\leq v-\rho^{'}.$$
For $a^{'}\rightarrow a$ and $b^{'}\rightarrow b$ we get by the
continuity of $\rho_{.}$ :
 $$\limsup_{t\rightarrow\infty}\frac{1}{t}\log\sharp G  (t)\leq v-\rho.$$
\end{proof}

\bigskip

Now we will give an other corollary, using the same method as that
of Bertoin and Gnedin in \cite{bergne}. We encode the configuration
$J^{\dag}  (t)=\{|J_{i}^{\dag}  (t)|\}$ of the lengths of good
intervals into the random measure
 $$\sigma_{t}:=\frac{e^{\rho t}}{h  (0)} \sum_{i\in\mathbb{N}}
h\left  (vt+\log|J_{i}^{\dag}  ( t)|\right)|J_{i}^{\dag}
 (t)|\delta_{\log  (1/a)+vt+\log  |J_{i}^{\dag}  (t)|}
$$
which has total mass $M_{t}$.

 The associated mean measure
$\sigma_{t}^{*}$ is defined by the formula
$$\int_{0}^{\infty}f  (x)\sigma_{t}^{*}  (dx)=\mathbb{E}\left (\int_{0}^{\infty}f  (x)\sigma_{t}  (dx)\right)$$
which is required to hold for all compactly supported continuous
functions $f$. Since $M_{t}$ is a martingale, $\sigma_{t}^{*}$ is a
probability measure. More precisely the next proposition establishes
the convergence of the mean measure $\sigma_{t}^{*} $, and then of
$\sigma_{t}$ itself.

\begin{proposition}\label{theorememesureempirique}
In the notation of the previous sections, with the assumptions~
 (\ref{nuinf}), and $v>\rho $ we get:
\begin{enumerate}
\item The measures $\sigma_{t}^{*}$ converge weakly, as
$t\rightarrow\infty$, to the probability measure
$$\varrho  (dy):=ch  (y+\log a)h  (\log  (b)-y)dy$$
where $c>0$ is the constant that appears in Theorem~\ref{theo11}.5.
\item For any bounded continuous $f$

\begin{equation}
L^{2}-\lim_{t\rightarrow\infty}\int_{0}^{\infty}f  (x)\sigma_{t}  (
dx)=M_{\infty}\int_{0}^{\infty} f  (x)\varrho  (dx).
\end{equation}
\end{enumerate}
\end{proposition}

\begin{proof}
\begin{enumerate}
\item Firstly we prove the convergence of the mean measures $\sigma_{t}^{*}\rightarrow\varrho$. Let $f$ be a bounded continuous function.
By definition we get:

$\int_{0}^{\infty}f  (y)\sigma_{t}^{*}  (dy)$

$$\begin{disarray}{rl}
=&\mathbb{E}\left  (\int_{0}^{1}f  (\log  (1/a)+vt+\log
|I_{x}^{\dag}
 (t)|)\frac{e^{\rho t}}{h  (0)} h\left  ( vt+\log|I_{x}^{\dag}  (
t)|\right)\mathbf{1}_{\{I_{x}^{\dag}  (t)\in G  (t)
\}}dx\right)\\&\\=&\mathbf{E}_{\log  (1/a)}\left  (f  (Y_{t})e^{\rho
t}\frac{h  (Y_{t}+\log a)}{h  (0)}\mathbf{1}_{\{t<T\}}\right),
\end{disarray}$$
with the definition of $Y_{t}$. Thus by the definition of $P_{t}$ in
Theorem~\ref{theo11}, we get
$$\int_{0}^{\infty}f  (y)\sigma_{t}^{*}  (dy)=\int_{0}^{\log   (b/a)}f  (y)\frac{h  (y+\log a)}{h  (0)}e^{\rho t}P_{t}  (\log  (1/a),dy).$$
By Theorem~\ref{theo11}.2, we get
$$\int_{0}^{\infty}f  (y)\sigma_{t}^{*}  (dy)\underset{t\rightarrow\infty}{\sim}c\int_{0}^{\log  (b/a)}f  (y)h  (y+\log a)h  (\log  (b)-y)dy.$$
Therefore the measure $\sigma_{t}^{*}$ converge weakly to the
probability measure $\varrho$.

\item Now we show that the scaled empirical measures induced by
$J  (t)$ converge in the $L^{2}$-sense to the random measure
$M_{\infty}\varrho$.

Let $f_{1}$ and $f_{2}$ be two continuous functions  bounded from
above by 1, and
$$\begin{disarray}{rcl}
S_{t}&=&\sum_{i,j}f_{1}  (\log  (1/a)+vt+\log |J_{i}^{\dag}   (
t)|)\frac{e^{\rho t}}{h  (0)} h\left  (vt+\log|J_{i}^{\dag}
 (t)|\right)|J_{i}^{\dag}  (t)|\
\\&&\\&&\times f_{2}  (\log  (1/a)+vt+\log |J_{j}^{\dag}  (t)|)\frac{e^{\rho
t}}{h  (0)} h\left  (vt+\log|J_{j}^{\dag}   (
t)|\right)|J_{j}^{\dag}
 (t)|.
\end{disarray}$$
We need to show that

\begin{eqnarray}
&&\mathbb{E}\left  (S_{t}\right) \rightarrow\left   (
\int_{0}^{\infty}f_{1}  (x)\varrho  (dx)\right)\left   (
\int_{0}^{\infty}f_{2}  (x)\varrho  (dx)\right)\mathbb{E}\left   (
M_{\infty}^{2}\right) \label{convergencel2}
\end{eqnarray}
for  $f_{1}$ and $f_{2}$ positive and bounded from above by 1.
Indeed, suppose~ (\ref{convergencel2}) is shown. Denote
$$A_{t}=\sum_{j}f_{1}  (\log  (1/a)+vt+\log |J_{i}^{\dag}  (t)|)\frac{e^{\rho
t}}{h  (0)} h\left  (vt+\log|J_{j}^{\dag}   (
t)|\right)|J_{j}^{\dag}
 (t)|.$$ Take $f_{2}=1$ to conclude from~ (\ref{convergencel2}) that
$$\lim_{t\rightarrow\infty}\mathbb{E}  (A_{t}M_{t})=\int_{0}^{\infty}f_{1}  (x)\varrho  (dx)\mathbb{E}\left  (M_{\infty}^{2}\right).$$
Similarly, by setting $f_{1}=f_{2}$ we get
$$\lim_{t\rightarrow\infty}\mathbb{E}\left  (A_{t}^{2}\right)=\left  (\int_{0}^{\infty}
f_{1}  (x)\varrho  (dx)\right)^{2}\mathbb{E}\left   (
M_{\infty}^{2}\right).$$ Recalling that $\mathbb{E}   (
M_{t}^{2})\rightarrow\mathbb{E}  (M_{\infty}^{2})$ and combining the
above we get the desired
$$\lim_{t\rightarrow\infty}\mathbb{E}\left[\left  (A_{t}-M_{t}\int_{0}^{\infty}f_{1}  (x)\varrho  (dx)\right)^{2}\right]=0.$$

To prove~ (\ref{convergencel2}) let us replace $t$ by $t+s$ and
condition on $J^{\dag}=  (|J_{i}^{\dag}  (s)|)_{i\in\mathbb{N}}$. We
have two cases: write $i\sim_{s} j$ for the case where at time $t+s$
two coexisting intervals $J_{i}^{\dag} (t+s)$ and $J_{j}^{\dag}
(t+s)$ stem from the same interval at time $s$, and $i\nsim_{s} j$
for the case these intervals  are not included into the same
interval component at time $s$. Therefore, with the notation
$$\begin{disarray}{rcl}S_{t+s}^{  (1)}:=\mathbb{E}\left
 (\sum_{i\sim_{s} j}\ S_{t+s}\ |\ J^{\dag}  (s)\right)& \hbox{and} &
S_{t+s}^{  (2)}:=\mathbb{E}\left  (\sum_{i\nsim_{s} j}\ S_{t+s}\ |\
J ^{\dag} (s)\right)\end{disarray}$$ we get:

$$S_{t+s}^{  (1)}+S_{t+s}^{  (2)}=\mathbb{E}  (S_{t+s}|J^{\dag}  (s)).$$

For the studies of $S_{t+s}^{  (1)}$ we use the homogeneous property
of the fragmentation and the notation $I_{0}=(0,\log  ( b/a))$, and
get

$|S_{t+s}^{  (1)}|$
$$\begin{disarray}{rl}
 \leq& \sum_{i}|J_{i}^{\dag}  (s)|^{2} e^{2\rho
s} \mathbb{E}\left  (\sum_{j} |J_{j}^{\dag}  (t)|e^{\rho
t}\right)^{2} \sup_{x\in I_{0}}\left  (\frac{h  (x+\log a)}{h  (
0)}\right)^{2}\sup_{x\in I_{0}}|f_{1}  (x)|\sup_{x\in I_{0}}|f_{2}
  (x)|\\&\\ \leq & be^{  (\rho-v)s} C_{6},\end{disarray}$$ with
\[C_{6}:=\sum_{i}|J_{i}^{\dag}  (s)|^{} e^{\rho t} \mathbb{E}\left
  (\sum_{j} |J_{j}^{\dag}  (t)|e^{\rho t}\right)^{2} \sup_{x\in
I_{0}}h  (x+\log  (a))^{2}\sup_{x\in I_{0}}|f_{1}  (x)|\sup_{x\in
I_{0}}|f_{2}  (x)|/h  (0)^{2}\]
which is finite because
$$\mathbb{E}\left  (\sum_{j} |J_{j}^{\dag}  (t)|e^{\rho
t}\right)=\mathbf{E}\left  (\mathbf{1}_{\{t<T \}}e^{\rho
t}\right)<\infty.$$ Thus $S_{t+s}^{  (1)}\rightarrow 0$ as
$s\rightarrow\infty$ uniformly in $t$.

Now we look at $S_{t+s}^{  (2)}$. We introduce the notation
$y_{k}=|J_{k}^{\dag}  (s)|$. Write $i\searrow k$ if the length
$|J_{i}^{\dag}  (t+s)|$ stems from $y_{k}$. By independence, the
intervals which are included in the interval with length $y_{k}$ and
those which are included in the interval with length  $y_{l}$ evolve
independently, thus gathering the lengths $|J_{i}^{\dag}   ( t+s)|$
by the ancestors at time $s$ yields
$$ S_{t+s}^{  (2)}= \sum_{k\neq l}\left  (\mathbb{E}\sum_{i\searrow k}\
...\right)\left  (\mathbb{E}\sum_{j\searrow l}\ ...\right).$$

On the other hand, by self-similarity and convergence of the mean
measures
$$\begin{disarray}{l}\mathbb{E}\left  (\sum_{i\searrow k}e^{\rho s}
\frac{h\left  (vs+\log\left  (y_{k}/a\right)\right)}{h  (0)}y_{k}\
f_{1}  (vt+\log  (|J_{i}^{\dag}  (t)|)+vs+\log  (y_{k}/a))\right.\\
\\ \left.\left.e^{\rho t}
\frac{h\left  (vt+\log\left  (|J_{i}^{\dag}  (t)|\right)+vs+\log  (
y_{k}/a)\right)}{h\left  (vs+\log\left   (
y_{k}/a\right)\right)}|J_{i}^{\dag}  (t)|\ \right|\ J^{\dag}   (
s)\right)
\\ \\\underset{t\rightarrow\infty}{\rightarrow}\ \ e^{\rho s}
\frac{h\left  (vs+\log\left  (y_{k}/a\right)\right)}{h  (
0)}y_{k}\left   (\int_{0}^{\infty}f_{1}  (x)\varrho  (dx)\right),
\end{disarray}$$
and
$$\begin{disarray}{l}\mathbb{E}\left  (\sum_{j\searrow l}e^{\rho s}
\frac{h\left  (vs+\log\left  (y_{l}/a\right)\right)}{h  (0)}y_{l}\
f_{2}   (vt+\log  (|J_{j}^{\dag}  (t)|)+vs+\log  (y_{l}/a))
\right.\\
\\ \left.\left.e^{\rho t}
\frac{h\left  (vt+\log\left  (|J_{j}^{\dag}  (t)|\right)+vs+\log  (
y_{l}/a)\right)}{h\left  (vs+\log\left   (
y_{l}/a\right)\right)}|J_{j}^{\dag}  (t)|\ \right|\ J^{\dag}   (
s)\right)
\\ \\\underset{t\rightarrow\infty}{\rightarrow}\ \ e^{\rho s}
\frac{h\left  (vs+\log\left  (y_{l}/a\right)\right)}{h  (
0)}y_{l}\left   (\int_{0}^{\infty}f_{2}  (x)\varrho  (dx)\right).
\end{disarray}$$
Therefore by dominated convergence
$$\begin{disarray}{rcl}\mathbb{E}\left  (S_{t+s}^{  (2)}\right)&
\underset{s\rightarrow\infty}{\sim}&\left.\left   (
\int_{0}^{\infty}f_{1}  (x)\varrho  (dx)\right)\left   (
\int_{0}^{\infty}f_{2}  (x)\varrho  (dx)\right) \mathbb{E} \right
  (\sum_{k\neq l} \frac{ e^{\rho s}}{h  (0)} |J_{k}^{\dag}  (s)|\\
&&\\&&\left. h\left  (vs+\log  (|J_{k}^{\dag}  (s)|/a)\right)
\frac{e^{\rho s}}{h  (0)} h\left  (vs+\log  (|J_{l}^{\dag}   (
s)|/a)\right)|J_{l}^{\dag}  (s)|\right).
\end{disarray}$$
Moreover with $C_{7}:= b \sup_{x\in I_{0}}|h  (x+\log a)|^{2}/h   (
0)^{2}$, we get
$$  \mathbb{E}\left  (\sum_{k}e^{2\rho s}
\frac{h\left  (vs+\log\left  (|J_{k}^{\dag}
 (s)|/a\right)\right)^{2}}{h  (0)^{2}}|J_{k}^{\dag}  (s)|^{2}\right)
\leq C_{7}\mathbb{E}\left  (\sum_{k}e^{\rho s} |J_{k}^{\dag}   (
s)|\right)e^{  (\rho-v)s}
$$ which goes to 0 when  $s\rightarrow\infty$, as a consequence
$$\begin{disarray}{l}
\mathbb{E} \left  (\sum_{k\neq l}    e^{\rho s} \frac{h\left   (
vs+\log\left  (|J_{k}^{\dag}  (s)|/a\right)\right)}{h   (
0)}|J_{k}^{\dag}  (s)| e^{\rho s}
\frac{h\left  (vs+\log\left  (|J_{l}^{\dag}  (s)|/a\right)\right)}{h  (0)}|J_{l}^{\dag}  (s)|\right)\\
\underset{s\rightarrow\infty}{\sim} \mathbb{E} \left   (
M_{s}^{2}\right).
\end{disarray}$$

\end{enumerate}

\end{proof}

\section{The Hausdorff dimension.}\label{sechausdorffdim}

In this section we use the notation and definitions of the previous
sections. We recall that $\rho=\rho_{\log  (b/a)}$, where $\rho_{.}$
is define in  (\ref{rho}).  Let $dim$ be the Hausdorff dimension.
The aim of this section would be to proof the main theorem:

\begin{theorem}\label{theoprincipal}\textbf{: Multifractal spectrum.} Assume~ (\ref{nuinf}):

$\bullet$ if $\rho > v $ holds, then:
$$G_{  (v,a,b)}=\emptyset\ \ \ a.s.$$

$\bullet$ if $\rho < v $ holds, then: \begin{equation} dim  (G_{
(v,a,b)})=1-\rho/v\ \ \ a.s. \label{dimG}
\end{equation}
\end{theorem}

\begin{rk}
\begin{enumerate}
\item
  Berestycki in \cite{be}  has computed the Hausdorff
dimension of the set $$G_{v}=\left\{x\in  (0,1) \ | \
\lim_{t\rightarrow\infty} \frac{1}{t}\log|I_{x}   (
t)|=-v\right\}.$$ He found that for $v\in  (v_{min},v_{max})$, $dim
 ( G_{v})=C  (v)/v$  (with $C  (v)$ defined at the beginning of
section~\ref{limittheorem}). In  Remark \ref{rk sur c  (v)} we have
shown that for all $v\in  (\max  (v_{min},\rho),v_{max})$ we have $C
(v)\geq v-\rho$ and we can notice that the inequality is strict for
$\rho >v_{min}$. As a consequence the set $G_{  (v,a,b)}$ has a
Hausdorff dimension smaller than that of $G_{v}$, and also  smaller
than that one could have infer from  equality  (\ref{eqberrou}).

\item In the case $v>v_{typ}$, we have $Y_{t}/
t\underset{t\rightarrow\infty}{\rightarrow} v-v_{typ}>0\ \ \ a.s.$
and \[\mathbf{P}_{\log  (1/ a)} (\inf\{t:\ Y_{t}\leq 0\}=\infty
)>0.\] Thus  $W^{ (-q)} (\infty )=0$ for all $q\geq0$ and then
$\lim_{\beta\rightarrow \infty}\rho_{\beta}=0$. Moreover using the
fact that, $\lim_{\beta\rightarrow 0}\rho_{\beta}=\infty$  and
$\rho_{.}$ is decreasing, we get that for all $v>v_{typ}$, there
exist $a$ and $b$ such that $\rho_{\log  (b/a)}<v$ and thus the fact
that the set of good intervals is not empty.
\end{enumerate}
\end{rk}

 The proof of this theorem use the non-asymptotic set $\Lambda_{  (v,a,b)}$. In  particular the key of the proof is the next proposition:

\begin{proposition}\label{theo6} Assume~ (\ref{nuinf}) and $0<a<b<1$:

$\bullet$ if $\rho > v $ holds, then:
$$\Lambda_{  (v,a,b)}=\emptyset\ \ \ a.s.$$

$\bullet$ if $\rho < v $ holds, then:  $\mathbb{P} (\Lambda_{
(v,a,b)}\neq \emptyset)>0$, and conditionally on $\Lambda_{
(v,a,b)}\neq \emptyset$,
\begin{equation}
dim  (\Lambda_{  (v,a,b)})=1-\rho/v. \label{dim}
\end{equation}
\end{proposition}

\medskip

\begin{proof}
\begin{enumerate}
\item Let $v>0$ and $a$ and $b$ such that $v<\rho$. We define
$$N  (t):=\sharp G  (t),$$with   $G  (t)$ defined in~ (\ref{good}).
 We remark that $$N  (t)=\int_{0}^{1}
\frac{1}{|I_{x}  (t)|}\mathbf{1}_{\{I_{x}  (t)\in G  (t)\}}  (x)dx
.$$

and in particular

 $$\begin{disarray}{rcl}
\mathbb{E}  (N  (t))&=&\mathbb{E}\left  (\int_{0}^{1}
\frac{1}{|I_{x}
 (t)|}\mathbf{1}_{\{I_{x}  (t)\in G  (t)\}}  (x)dx\right) .
\end{disarray}$$

Additionally by~ (\ref{xi}), we get

$$\begin{disarray}{rcl}
\mathbb{E}  (N  (t))&=& e^{vt}\ \mathbf{E}\left  (e^{\xi (t)-vt}\
\mathbf{1}_{\{vs-\xi (s)-\log a\in  (0,\log  (b /a))\ \forall\
s\leq t\}}\right).\\
 \end{disarray}$$

With the notation $Y_{t}=vt-\xi (t)$ and $P_{t}$ defined in
Theorem~\ref{theo11} we rewrite the previous equality as:

$$\begin{disarray}{rcl}
\mathbb{E}  (N  (t))&=&e^{vt}\
\mathbf{E}_{\log  (1/a)}\left  (e^{-Y_{t}-\log a}\ \mathbf{1}_{\{t<T\}}\right)\\&&\\&=&\frac{1}{a}e^{  (v-\rho)t}\ \int_{0}^{\log  (b/a)} e^{-y+\rho t} P_{t}  (\log  (1/a), dy).\\
 \end{disarray}$$

By Theorem~\ref{theo11}.2  we get

$$\begin{disarray}{rcl}
\mathbb{E}  (N  (t))&\underset{t\rightarrow\infty}{\sim} &
\frac{1}{a}e^{  (v-\rho)t}\ c\ h  (0) \ \int_{0}^{\log  (b/a)}
e^{-y}\ \ \ \Pi  (dy),
 \end{disarray}$$ with $\Pi$ defined in Theorem~\ref{theo11}.

Finally as the function  $y\mapsto e^{-y} \ h  (\log  (b)-y)$ is
continuous, the integral above is a finite constant. Thus if
$\rho>v$ then $\lim_{t\rightarrow\infty} \mathbb{E}  (N  (t))=0$,
from which one concludes that $\lim_{t\rightarrow\infty} N  (t)=0$ ,
i.e. $\Lambda_{  (v,a,b)}=\emptyset$ $a.s.$

\item Now we deal with the case where $a$ and $b$ are such that $v>\rho$. We work conditionally on $\Lambda_{ (v,a,b)}\neq \emptyset$  (or, equivalently, on the event $\zeta=\infty$, which has a positive probability by Theorem \ref{theo2}).

$\bullet$ Firstly, in order to prove the lower bound of the
Hausdorff dimension of $\Lambda_{  (v,a,b)}$, we will use the same
method as Berestycki in \cite{be} . We will divide this proof into
three steps. Each step will begin with a star $  (\star)$.  In the
first step we will construct a subset
$\underset{n\in\mathbb{N}}{\cap}\mathbb{G}_{\delta}  (n)$ of
$\Lambda_{  (v,a,b)}$, which will be defined latter on (see
(\ref{Gdelatn})). In the second we shall obtain a lower bound of the
Hausdorff dimension of this subset. In order to do that we will
construct an increasing process indexed by $t\in (0,1)$, which only
increases on $\underset{n\in\mathbb{N}}{\cap}\mathbb{G}_{\delta}
(n)$, and which is Hölder continuous.  In the last step we will
conclude.

$\star$ As in \cite{be}  for $ \delta>0$ we define for all $n\in
\mathbb{N}$, $H_{\delta}  (n)$ as a multi-type branching process
with each particle corresponding to a segment of $G  (\delta n)$ and
$$G_{\delta}  (n):=\underset{I\in H_{\delta}  (n)}{\cup} I,$$
with   $G  (t)$ defined in~ (\ref{good})  (i.e. $G_{\delta}
 (n)=G (\delta n))$.

We notice that the family $  (G_{\delta}   ( n))_{n\in\mathbb{N}}$
is nested and that $\underset{n\in\mathbb{N}}{\cap}G_{\delta}
 (n)=\Lambda_{  ( v,a,b)}$.

Let $\epsilon>0$, and fix $\epsilon^{'}>0$ and $\eta>0$ such that
$\eta <\min  (\epsilon, v-\rho ).$ By Proposition~\ref{theo4}, for
this $\epsilon^{'}>0$ and $\eta>0$, we may find $t_{0}>\max ( (1+|
\log  (1-\epsilon^{'})|)/ (\epsilon -\eta ),\log  (2)/ (v-\rho -\eta
))$ such that for all $t>t_{0}$:
$$\mathbb{P}  (|t^{-1}\log  (\sharp
G  (t))-  (v-\rho)|>\eta|\zeta=\infty)<\epsilon^{'}.$$ For each
$t>0$, we consider a variable $ \overset{\sim}{\chi}  (t)$ whose law
is given by$$ \mathbb{P}  (\overset{\sim}{\chi}
(t)=0)=\epsilon^{''},$$ and
$$
\mathbb{P}  (\overset{\sim}{\chi}  (t)=\lfloor e^{[   (
v-\rho)-\eta]t}\rfloor )=1-\epsilon^{''},$$  where  $\lfloor .
\rfloor$ is the  integer part
 and
$\epsilon^{''}:=\mathbb{P}  (|t^{-1}\log  (\sharp G  (t))-   (
v-\rho)|>\eta|\zeta=\infty)<\epsilon^{'}$. Moreover by using that
for all $x\geq 2$: $\log  ( x)-1\leq \log  (\lfloor x\rfloor )$, we
notice that
$$|t^{-1}\log  (\mathbb{E}  (\overset{\sim}{\chi}  (t)))-  (v-\rho)|\leq
\eta+t^{-1}  (|\log  (1-\epsilon^{'})|+1).$$ Plainly
$\overset{\sim}{\chi}  (t)$ is stochastically dominated by $\sharp G
(t)$. Exactly as in \cite{be}  we can construct a true Galton-Watson
tree $ \mathbb{H}$ by thinning $H_{\delta}$ where $\delta>t_{0}$.
More precisely the offspring distribution of $\mathbb{H}$ is given
by the law of $\overset{\sim}{\chi}   ( \delta)$. Let $m:=
\mathbb{E}  (\overset{\sim}{\chi}  (\delta))$ be the expectation of
the number of children of a particle. Therefore, we get
\begin{enumerate}
\item
\begin{equation}
|\delta^{-1}\log m -  (v-\rho)|<\epsilon. \label{14}
\end{equation}
\item The family   \begin{equation}
  (\mathbb{G}  (n):=\underset{I\in \mathbb{H}  (n)}{\cup}
 I)_{n\in \mathbb{N}} \label{Gdelatn}
\end{equation} is nested. The  $
 \mathbb{G}  (n)$ is the union of the interval of the $n$ generation of $\mathbb{H}$.
\item
$\underset{n\in\mathbb{N}}{\cap}\mathbb{G}  (n)\subseteq\Lambda_{
  (v,a,b)}.$

This last point makes sense because we work conditionally on
$\zeta=\infty$.
\end{enumerate}

$\star$ We fix $\epsilon> 0$. We choose $\delta>t_{0}$ as shown
above and consider the tree $ \mathbb{H}$. We define $Z (n)$ as the
number of nodes of  $\mathbb{H}$ at height $n$. By the theory of
Galton-Watson processes, as we are working conditionally on the
event $\Lambda_{ (v,a,b)}\neq \emptyset$, we have that almost surely
$$m^{-n}Z  (n)\rightarrow \mathcal{W}>0.$$
Let $\sigma$ be a node of our tree  (thus it is also a subinterval
of $ (0,1)$). Fix an interval $I\subset  (0,1)$ and introduce
$$\begin{disarray}{rcl}
 \mathbb{ H}_{I}  (n)&:=&\{\sigma\in\mathbb{H}  (n), \sigma\cap I\neq
 \emptyset\},\\&&\\Z_{I}  (n)&:=&\sharp\ \mathbb{H}_{I}  (n).
\end{disarray}$$
Define
\[x\rightarrow L_{x}:=\lim_{n}m^{-n}
Z_{ (0,x)}  (n),\ \  x\in\  (0,1).\] We will now state a lemma that
we will use to conclude:

\begin{lemma}\label{lemma3}
For each $\epsilon>0$, \begin{enumerate}
\item There exists a version $\overset{\sim}{L}$ of $ (L_{x})_{x\in[0,1]}$ which is
Hölder continuous of order $\alpha$ for any $\alpha<  1
-\rho/v-\epsilon$ for every $\epsilon>0$.
\item The process $\overset{\sim}{L}$ only grows on the set
$\underset{n\in\mathbb{N}}{\cap}\mathbb{G}  (n)$.
\end{enumerate}
\end{lemma}

\begin{proof}[Proof of Lemma \ref{lemma3}.]
\begin{enumerate}\item Exactly as in \cite{be} , we show the first point by
verifying Kolmogorov's criterium  (see \cite{revyor} Theorem 2.1
p.26). Let $W (\sigma )$ be the ``renormalized weight'' of the tree
rooted at $\sigma$, i.e.,
$$W (\sigma ):=\lim_{n\rightarrow \infty}m^{-n}\ \sharp\{\sigma^{'}\in\mathbb{H}  (|\sigma |+ n), \sigma^{'}\subset \sigma \},$$ where $|\sigma |$ is the generation of $\sigma$.

By the definition of $L$ we have for all $x>y\in  (0,1)$:
\[|L_{x}-L_{y}|=\lim_{n\rightarrow\infty} m^{-n}
Z_{ (x,y)}  (n),\ \  x\in\  (0,1).\]

For any $J$ open subinterval of $ (0,1)$, let
$$\eta (J):=\sup\{n\in\mathbb{N}:\ e^{-v\delta n}\geq |J|\}=\lfloor -\log (|J|)/ v \delta \rfloor.$$

For all $x,y$ such that $x<y$ by the definition of $L_{.}$, we get:
$$\begin{disarray}{l}
|L_{x}-L_{y}|\\=\lim_{n} m^{-\eta ( (x,y))} m^{-n+\eta ( (x,y))}\sum_{\sigma \in \mathbb{H}_{ (x,y)} (\eta  ( (x,y)))}  \sharp\{\sigma^{'}\in\mathbb{H}  (|\sigma |+n- \eta ( (x,y))), \sigma^{'}\subset \sigma \}\\
\leq   m^{-\eta ( (x,y))} \sum_{\sigma \in \mathbb{H}_{ (x,y)} (\eta
 ( (x,y)))}W (\sigma), \end{disarray}$$ and by the definition of
$\eta (.)$:
$$\begin{disarray}{rcl}
|L_{x}-L_{y}|&\leq &  e^{\log m (\frac{1}{v\delta }\log  (y-x)+1)} \sum_{\sigma \in \mathbb{H}_{ (x,y)} (\eta  ( (x,y)))}W (\sigma )\\
&\leq & m |x-y|^{1-\epsilon-\rho/ v} \sum_{\sigma \in \mathbb{H}_{
(x,y)} (\eta  ( (x,y)))}W (\sigma ), \end{disarray}$$ by using
(\ref{14}). Moreover by the definition of good intervals, we have
that for each $n$ the sizes of  intervals in $\mathbb{H} (n)$ have a
lower bound given by $ae^{-v\delta n}$, so $a|J|e^{-v\delta}$ is a
lower bound for the sizes of the intervals of $\mathbb{H} (\eta
(J))$, and thus $Z_{J} (\eta (J))\leq e^{v\delta}/a$. Therefore for
all $\gamma> 1$ and all $J\subset  (0,1)$ we have:
$$\begin{disarray}{rcl}\mathbb{E}\left ( (\sum_{\sigma \in \mathbb{H}_{J} (\eta  (J))}W (\sigma ))^{\gamma
}\right)&\leq&
 \mathbb{E} ( (W_{1}+...+W_{\lfloor e^{v\delta}/a \rfloor+1})^{\gamma })\\&\leq& \mathbb{E} ( (W_{1}+...+W_{\eta  (J)+2})^{\gamma })<\infty,\end{disarray}$$
 where the $W_{i}$ are i.i.d. with the same law as $W$. The finiteness comes from the existence of finite moments of all orders for $W$  (see for example
 Theorem 3.4 p. 479 of Harris \cite{har1}).

\item The second point is clear by the choice of $L$.
\end{enumerate}
\end{proof}

$\star$ To prove that $dim  \left (
\underset{n\in\mathbb{N}}{\cap}\mathbb{G}  (n)\right)\geq  1
-\rho/v-\epsilon$, it is enough to show that
\begin{equation}
\sum_{i} \  diam  (U_{i})^{ 1 -\rho/v-\epsilon}>0
\label{sommeaverifierpourborneinf}
\end{equation}
for any cover $\{U_{i}\}$ of
$\underset{n\in\mathbb{N}}{\cap}\mathbb{G}  (n)$, where $diam   (
U_{i})$ is the diameter of $U_{i}$. Clearly, it is enough to assume
that the $\{U_{i}\}$ are intervals, and by expanding them slightly
and using the compactness of the closure of
$\underset{n\in\mathbb{N}}{\cap}\mathbb{G}  (n)$, we only need to
check~ (\ref{sommeaverifierpourborneinf}) if $\{U_{i}\}$ is a finite
collection of open subintervals of $[0,1]$.

Let $\cup_{i=0}^{N}  (l_{i},r_{i})$ be a cover of
$\underset{n\in\mathbb{N}}{\cap}\mathbb{G}  (n)$  (where the $
(l_{i},r_{i})$ are disjoints open intervals). Therefore
$$\sum_{i=1}^{N}|\overset{\sim}{L}_{r_{i}}-\overset{\sim}{L}_{l_{i}}|=\mathcal{W}.$$
Thus for all such covers with $\max_{i}  (r_{i}-l_{i})$ small enough
$$\mathcal{W}\leq k\sum_{i=0}^{N}  (r_{i}-l_{i})^{ 1
-\rho/v-\epsilon}$$ and hence
$$dim  (\Lambda_{  (v,a,b)})\geq dim  (\underset{n\in\mathbb{N}}{\cap}\mathbb{G}  (n))\geq
 1
-\rho/v-\epsilon .$$ To get the lower bound of the Hausdorff
dimension of $\Lambda_{  (v,a,b)}$,  we let $\epsilon$ tend to 0.

\medskip

$\bullet$  Secondly, the upper bound for~ (\ref{dim}) is an easy
corollary of the fact that the Hausdorff dimension is smaller than
the box-counting dimension  (see \cite{fal} p.36-43), using the
cover $\underset{n\geq N}{\cup}\ \ \underset{i\in\theta_{v,a,b}   (
n)}{\cup} J_{i}  (n)$, with $\theta_{v,a,b}   (
t)=\left\{i\in\mathbb{N}  \ |\ J_{i}  (t)\in G  (t)\right\}$  (with
$ G   (t)$ defined in Section~\ref{sec}).

\end{enumerate}
\end{proof}

Then we have the next corollary, which deals with the general case
for $a$ and $b$:

\begin{corollary}\label{cor1} For $t^{'}\geq 0$  set $$ \Lambda_{  (v,a,b)} (t^{'}):=\left\{x\in  (0,1)\ : \ a e^{-vt}<|I_{x}  (t)|<be^{-vt} \ \forall t\geq t^{'}\right\}.$$ Assume  (\ref{nuinf}), $0<a<b$ and $\rho<v$, then
 \[\mathbb{P} (\Lambda_{ (v,a,b)} (t^{'})\neq \emptyset)\underset{t^{'}\rightarrow \infty}{\rightarrow} 1,\]
 and
\[\mathbb{P}\left.\left (dim (\Lambda_{ (v,a,b)} (t^{'}))=1-\rho / v\ \ \right|\Lambda_{ (v,a,b)} (t^{'})\neq \emptyset\right)=1.\]
\end{corollary}

\begin{proof}

\begin{enumerate} \item The first part of the proof is a consequence of the homogeneity of the fragmentation and of Proposition \ref{theo6}.

\item Fix $\rho^{'}>\rho$.  As
 $\lim_{\beta \rightarrow 0} \rho_{\beta}=\infty $, and, by Theorem \ref{theo11}.5, the application $\beta\rightarrow \rho_{\beta}$ is continuous and strictly decreasing, therefore there exists $\beta_{0}\in  (1,b/a)$ such that $\rho^{'}=\rho_{\log  (\beta_{0})}$.
Let $\epsilon:= (\beta_{0}-1)/ (1+\beta_{0})$, $a^{'}:=1-\epsilon$,
$b^{'}:=1+\epsilon$, $x_{0}:= (\beta_{0}+1) (a+b/\beta_{0})/4$
(notice that $x_{0}\in  (a,b)$) and
$$p_{0}:=\mathbb{P} (dim (\Lambda_{ (v,a^{'},b^{'})})\geq1-\rho_{\log  (b^{'}/a^{'})}/ v). $$
By Proposition \ref{theo6}, we get that $p_{0}>0$. We notice that by
the choice of $a^{'}$ and of $b^{'}$, we have  $\rho_{\log
(b^{'}/a^{'})}=\rho_{\log  (\beta_{0})}=\rho^{'}$.

Let $I$ be an interval of $ (0,1)$. The law of the homogeneous
interval fragmentation started at $I$ will be denoted by
$\mathbb{P}_{I}$. We remark that $\mathbb{P}_{I}  (dim (\Lambda_{
(v,a,b)})\geq1-\rho^{'}/ v) $ only depends on the length of $I$.
Thus we define $$g_{a,b}  (x):=\mathbb{P}_{I}   (dim (\Lambda_{
(v,a,b)})\geq1-\rho^{'}/ v),$$ where $I$ is an interval such that
$|I|=x$.

Let $x\in  (x_{0} a^{'},x_{0} b^{'})$. We remark that by the choice
of $x_{0}$ and as $1<\beta_{0}<b/a$ we have that $  (x_{0}
a^{'},x_{0} b^{'})\subset  (a,b)$ and thus
$$g_{a,b}  (x)\geq g_{x_{0} a^{'},x_{0} b^{'}}  (x).$$
Moreover by the scaling property of the fragmentation we get that $$
g_{x_{0} a^{'},x_{0} b^{'}}  (x)=\mathbb{P} (dim (\Lambda_{
(v,a^{'}/x,b^{'}/x)})\geq1-\rho_{\log  ( (b^{'}/x)/ (a^{'}/x))}/ v)=
p_{0}$$ Therefore
\begin{equation}\inf_{x\in  (x_{0} a^{'},x_{0} b^{'})}g_{a,b}  (x)\geq  p_{0}.\label{eqp}\end{equation}

 Let $$B (t)= \{i :\ \ x_{0} a^{'}< e^{vt }|J_{i } (t)|<x_{0}b^{'}\} \ \ ,\    n_{t}=\sharp B (t),$$ where $ (J_{1},J_{2},...)$ is the interval decomposition of $F (t)$.

 Fix $t^{'}\geq 0$. By applying the Markov property at time $t^{'}$  we get that
$$\begin{disarray}{l}\mathbb{P} (dim (\Lambda_{ (v,a,b)} (t^{'}))< 1-\rho^{'}/ v))\\\leq  \mathbb{E}\left (\prod_{i\in B (t^{'})}\mathbb{P}_{J_{i} (t^{'})} (dim (\Lambda_{ (v,x_{0}a^{'},x_{0}b^{'})})<1-\rho^{'}/v)\right)\\
\leq \mathbb{E} ( (1-p_{0})^{n_{t^{'}}}),
 \end{disarray}$$
 by using  (\ref{eqp}).
Therefore as $p_{0}>0$,
$n_{t^{'}}\underset{t^{'}\rightarrow\infty}{\rightarrow}\infty$
 (see  (\ref{eqberrou})) and with the first part of the proof we can
conclude.
\end{enumerate}
\end{proof}

\medskip

Now we are able to proof our main result:

\medskip

\begin{proof}[Proof of Theorem \ref{theoprincipal}.]

Observe that for all $n\in \mathbb{N}$, we have
\begin{equation}\Lambda_{  (v,a,b)} (n)\subset G_{  (v,a,b)}\subset \underset{\epsilon>0}{\cap}\underset{m\in\mathbb{N}}{\cup}
\Lambda_{  (v,a-\epsilon,b+\epsilon)} (m). \label{eqensembliste}
\end{equation}
We can notice that the second inclusion is actually an equality.

$\bullet$ First we consider the case where $\rho >v$. As the
application $\beta\rightarrow \rho_{\beta}$ is continuous and
strictly decreasing  (see Theorem \ref{theo11}.5),  there exists
$\epsilon_{0} > 0$ such that $v<\rho_{\log  ( (b+\epsilon_{0} )/
(a-\epsilon_{0} ))}<\rho$. Moreover by  (\ref{eqensembliste}) $$ G_{
(v,a,b)}\subset \underset{m\in\mathbb{N}}{\cup} \Lambda_{
(v,a-\epsilon_{0},b+\epsilon_{0})} (m),$$therefore thanks to
Proposition \ref{theo6} and the homogeneous property of the
fragmentation, we get the first part of the proof.

$\bullet$ Second we consider the case where $\rho <v$. Thanks the
second inclusion and the corollary \ref{cor1}, we get that: for all
$\epsilon\in  (0,a)$,
$$dim (G_{ (v,a,b)})\leq dim  (\underset{n\in\mathbb{N}}{\cup}
\Lambda_{  (v,a-\epsilon,b+\epsilon)} (n))=\max_{n} dim (\Lambda_{
(v,a-\epsilon,b+\epsilon)} (n))\leq 1-\rho_{\log
(\frac{b+\epsilon}{a-\epsilon})}/v.$$ Then by the continuity of
$\rho_{.}$  (see Theorem \ref{theo11}.5), we get the uper bound of
the Hausdorff dimension of $G_{ (v,a,b)}.$

The lower bound of the Hausdorff dimension is a consequence of the
first inclusion of  (\ref{eqensembliste}), as $dim  (\Lambda_{
(v,a,b)} (n))=1-\rho /v$ with a probability which  goes to 1 when
$n$ goes to infinity.

\end{proof}

\section{Appendix}

\subsection{A partition fragmentation.}
In this appendix we give a proof of Theorem \ref{theo2}.1.
(Section~\ref{sec}).

For this, we use the method of Bertoin and Rouault in \cite{berrou}
for fragmentation, which goes back to Lyons and al. \cite{lyoal} for
Galton-Watson processes, and  tools taken from the article of
Engl\"ander, Harris and Kyprianou \cite{enhaky}.

We first introduce the notations that we  need and we define what a
partition fragmentation $\Pi$ is. Let $ \mathcal{P}$ the space of
partition of $ \mathbb{N}$, and for every integer $k$, the block
$\{1,...,k\}$ is denoted by $[k]$. As in \cite{berrou}, we call
discrete point measure on the space
$\Omega:=\mathbb{R}_{+}\times\mathcal{P}\times\mathbb{N}$, any
measure :
$$w=\sum_{ (t,\pi,k)\in\mathcal{D}}^{\infty} \delta_{ (t,\pi,k)},$$
where $\mathcal{D}$ is a subset of
$\mathbb{R}_{+}\times\mathcal{P}\times\mathbb{N}$ such that
$$\forall t^{'}\geq 0\ \ \forall n\in\mathbb{N}\ \
\sharp\left\{ (t,\pi,k)\in\mathcal{D} \ |\ t\leq t^{'},
\pi_{|[n]}\neq ([n],\emptyset,\emptyset,...),k\leq n\right\}<\infty
$$
and for all $t\in\mathbb{R}$
\[w(\{t\}\times \mathcal{P}\times \mathbb{N})\in\{0,1\}.\]

Starting from an arbitrary discrete point measure $\omega$ on
$\mathbb{R}_{+}\times\mathcal{P}\times\mathbb{N}$, we will construct
a nested partition $\Pi= (\Pi (t),t\geq 0)$ (which means that for
all $t\geq t^{'}$ $\Pi (t)$ is a finer partition of $ \mathbb{N}$
than $\Pi (t^{'})$). We fix $n\in\mathbb{N}$, the assumption that
the point measure $\omega$ is discrete enables us to construct a
step path $ (\Pi (t,n),t\geq 0)$ with values in the space of
partitions of $[n]$, which only jumps at times $t$ at which the
fiber $\{t\}\times\mathcal{P}\times\mathbb{N}$ carries an atom of
$\omega$, say $ (t,\pi,k)$, such that $\pi_{|[n]}\neq
([n],\emptyset,\emptyset,...)$ and $k\leq n$. In that case, $\Pi
(t,n)$ is the partition obtained by replacing the $k-th$ block of
$\Pi (t-,n)$, denoted $\Pi_{k} (t-,n)$, by the restriction
$\pi_{|\Pi_{k} (t-,n)}$ of $\pi$ to this block, and leaving the
other blocks unchanged. Of course for all $t\geq 0$, $ (\Pi
(t,n),n\geq 0)$ is compatible (i.e. for every $n$, $\Pi (n,t)$ is a
partition of $[n]$ such that the restriction of $\Pi (n+1,t)$ to
$[n]$ coincide with $\Pi (n,t)$), as a consequence, there exists a
unique partition $\Pi (t)$, such that for all $n\geq 0$ we have $\Pi
(t)_{|[n]}=\Pi (t,n)$. With the terminology of \cite{ber2}, it is shown  in
\cite{berrou} that this process $\Pi$ is a
 (partition valued) homogeneous fragmentation.

One says that a block $B\subset \mathbb{N}$ has an asymptotic
frequency, if the limit
$$|B|:=\lim_{n\rightarrow\infty}n^{-1} card (B\cap [n])$$exists.
 When every block of some partition $\pi\in \mathcal{P}$
 has an asymptotic frequency, we write
 $|\pi|= (|\pi_{1}|,...)$ and then $|\pi|^{\downarrow}= (|\pi_{1}|^{\downarrow},...)\in\mathcal{S}^{\downarrow}$ for the
 decreasing rearrangement of the sequence $|\pi|$.
 In the case where some block of the partition $\pi$ does
 not have an asymptotic frequency, we decide that
 $|\pi|=|\pi|^{\downarrow}=\partial$, where $\partial$
 stands for some extra point added to $ \mathcal{S}^{\downarrow}$.
\textbf{We stress that the process of ranked asymptotic
frequencies $|\Pi|^{\downarrow}$ is a ranked fragmentation.}

Moreover, let $\nu$ be the dislocation measure associated to this
ranked fragmentation (see Subsection~\ref{pre}). According to
Theorem 2 in \cite{ber2}, there exists a unique measure $\mu$ on $
\mathcal{P}$ which is exchangeable (i.e. invariant by the action of
finite permutations on $ \mathcal{P}$), and such that $\nu$ is the
image of $\mu$ by the map  that associate the decreasing
rearrangement $|\pi|^{\downarrow}$ of the sequence of the asymptotic
frequencies of the blocks of $\pi$, to $\pi$. Thanks to
exchangeability, we get that for all measurable function $f:
[0,1]\rightarrow\mathbb{R}_{+}$  such that $f (0)=0$.
$$\int_{\mathcal{P}}f (|\pi_{1}|)\mu (d\pi)=\int_{\mathcal{S^{*}}}\sum_{i=1}^{\infty}s_{i}f (s_{i})\nu (ds).$$

 We denote the sigma-field generated by the
restriction to $[0,t]\times\mathcal{P}\times\mathbb{N}$ by $
\mathcal{G}_{0} (t)$. So $ (\mathcal{G}_{0} (t))_{t\geq 0}$ is a
filtration,  and the nested partitions $ (\Pi (t),t\geq 0)$ are $
(\mathcal{G}_{0} (t))_{t\geq 0}$-adapted. We define also  the
sigma-field  $ (\mathcal{F}_{0} (t))_{t\geq 0}$ generated by the
decreasing rearrangement $|\Pi (r)|^{\downarrow}$ of the sequence of
the asymptotic frequencies of the blocks of $\Pi (r)$ for $r\leq t$.
Of course $ (\mathcal{F}_{0} (t))_{t\geq 0}$ is a sub-filtration of
$ (\mathcal{G}_{0} (t))_{t\geq 0}$.

 Let
$ \mathcal{G}_{1} (t)$ the sigma-field generated by the
restriction of the discrete point measure $w$ to the fiber
$[0,t]\times\mathcal{P}\times\{1\}$. So $ (\mathcal{G}_{1} (t),t\geq
0)$ is a sub-filtration of $ (\mathcal{G}_{0} (t),t\geq 0)$, and the
first block of $\Pi$ is $ (\mathcal{G}_{1} (t),t\geq 0)$-measurable.
Let $\mathcal{D}_{1}\subseteq \mathbb{R}_{+}$ be the random set of
times $r\geq 0$ for which the discrete point measure has an atom on
the fiber $\{r\}\times\mathcal{P}\times\{1\}$, and for every
$r\in\mathcal{D}_{1}$, denote the second component of this atom by
$\pi (r)$.

 We define the probability measure $\mathbf{P}^{\updownarrow}$ as the $h$-transform of $\mathbf{P}$ based on
 the martingale $D_{t}$ (defined in
Theorem~\ref{theo11}  (3)):
\begin{equation}d\mathbf{P}_{x}^{\updownarrow} |_{\mathcal{E}_{t}}\
=D_{t} d\mathbf{P}_{x}
|_{\mathcal{E}_{t}}.\label{changementproba}\end{equation}

To simplify the notation, as in the section~\ref{sec} we define for
all $t\in\mathbb{R}$
 $h (t)=
W^{ (-\rho)} (t+\log (1/a))\mathbf{1}_{\{t\in (\log (a),\log
(b))\}}$. This function is well defined thanks to
Theorem~\ref{theo11}.

Let $P_{i} (t)$ the block of $\Pi (t)$ which contains $i$ at time
$t$. Similarly as in Section~\ref{sec}, for a homogeneous
fragmentation, we define the killed partition \[\Pi^{\dag}_{j}
(t)=\Pi_{j} (t)\mathbf{1}_{\{\exists i \in \mathbb{N}^{\ast}| \
\Pi_{j} (t)=P_{i} (t);\ \forall s\leq t\ \  |P_{i} (s)|\in
(ae^{-vs},be^{-vs})\}}.\]

 When we project
the martingale $D_{t}$ of~ (\ref{dt}) on the sub-filtration $
(\mathcal{G}_{0} (t))_{t\geq 0}$, we obtain an additive martingale
$$ \frac{e^{\rho
t}}{h(0)}\sum_{i=1}^{\infty} h (vt+\log (|\Pi_{i}^{\dag} (t)|))\
|\Pi_{i}^{\dag} (t)|\ .$$ As $|\Pi|$ is a ranked fragmentation with
dislocation measure $\nu$, this martingale is the same as this of
Section~\ref{sec}. From now on, we denote this martingale by $M_{t}$
too.

Observe that the projection (\ref{changementproba}) on the
sub-filtration $\mathcal{G}_{0} (t)$ give the identity:
\[d\mathbb{P}_{x}^{\updownarrow} |_{\mathcal{G}_{0}(t)}\
=M_{t} d\mathbb{P}_{x} |_{\mathcal{G}_{0}(t)}.\]

Like in  lemma 8  (ii) \cite{berrou}, with
 the probability measure $\mathbb{P}^{\updownarrow}$
we get:
\begin{lemma}\label{lm1}
Under $\mathbb{P}^{\updownarrow}$, the
 restriction of $w$ to
 $\mathbb{R}_{+}\times\mathcal{P}\times\{2,3,...\}$ has the same
 distribution as under $\mathbb{P}$ and is independent of the
 restriction to the fiber
 $\mathbb{R}_{+}\times\mathcal{P}\times\{1\}.$
\end{lemma}

It follows immediately from Theorem~\ref{theo11} that

\begin{rk}\label{lm2}
 For $x\in[0,\log (b/a)]$, let $F_{x} (t):=\mathbf{E}_{x}\left (\ e^{\rho t}
\mathbf{1}_{\{T>t\}\}}\right)$ for $t\in[0,\infty)$, then $F_{x}
(t)$ converges when $t\rightarrow\infty$ to a finite limit, and
$F_{x} (.):[0,\infty)\rightarrow[0,\infty)$ is c\`adl\`ag. In
particular we have \[\sup_{x\in[0,\log (b/a)]}\sup_{t\geq 0}|F_{x}
(t)|<\infty .\]
\end{rk}

\begin{rk}\label{rkmt} We have for all $t\geq 0$:
\[|M_{t}-M_{t-}|\leq e^{(\rho -v)t}\frac{b^{2}}{a h(0)}\sup_{x\in[\log a,\log
b]}h(x)\ \ \ \ a.s.\] If $v>\rho$,  there exists $0<C^{'}<\infty$
such that
 \[\sup_{t\geq 0}|M_{t}-M_{t-}|<C^{'}\  \  \text{a.s.}\]
\end{rk}

 Let
$$c_{t}:=\frac{e^{\rho t}}{h (0)}h\left (vt+\log (|\Pi_{1}^{\dag} (t)|)\right)
\ |\Pi_{1}^{\dag} (t)|$$ and
$$d_{t}:=\frac{e^{\rho t}}{h (0)} \sum_{i=2 }^{\infty}
h\left (vt+\log\left (|\Pi^{\dag}_{i} (t)|\right)\right)\
|\Pi^{\dag}_{i} (t)|.
$$

\bigskip

Now we have the background that we need to study

$$M_{t}=c_{t}+d_{t},$$
and we will show that $M$ is bounded in
$\mathrm{L}^{2}(\mathbb{P})$. In order to do that, as $\mathbb{E
}(M_{t}^{2})=\mathbb{E}^{\updownarrow}(M_{t})$, it is enough to
prove that
$$\lim_{t\rightarrow \infty}\mathbb{E}^{\updownarrow} (M_{t})<\infty\
 .$$

\subsection{The proof of Theorem
\ref{theo2}.1.}

$\bullet$ First we show that $\lim_{t\rightarrow \infty}
\mathbb{E}^{\updownarrow} (c_{t})=0 . $

With the subordinator $\xi (t):=-\log (|\Pi_{1} (t)|)$,   whose
Laplace exponent is $\kappa$ (exactly the same as this defined in
Subsection~\ref{subordinateur}), with the L\'evy Process
$Y_{t}=vt-\xi (t)+\log (1/a)$, and $T:=T_{\log (b/ a)}$ defined in~
(\ref{ta}) associated to this L\'evy Process, under
$\mathbf{P}_{\log (1/a)}^{\updownarrow}$ we get:
 $$\begin{disarray}{rcl}
c_{t}&=& \frac{e^{ (\rho-v) t}}{h (0)}W^{ (-\rho)} (Y_{t})\
e^{Y_{t}}\ \mathbf{1}_{\{t<T\}}.
\end{disarray}$$
As a consequence,
 $$\begin{disarray}{l}\mathbf{E}_{\log (1/a)}^{\updownarrow}\left(\frac{W^{ (-\rho)} (Y_{t})}{h(0)}\
e^{Y_{t}}\ \mathbf{1}_{\{t<T\}}\right)= \mathbf{E}_{\log
(1/a)}\left(\frac{W^{ (-\rho)}
 (Y_{t})^{2}}{h(0)^{2}}\ e^{Y_{t}}\ e^{\rho t}\mathbf{1}_{\{t<T\}}\right)\\\leq \sup_{x\in[\log a, \log b]}\left(h(x)\right)^{2}\frac{b}{ah(0)^{2}}\ F_{\log (1/a)}(t)\end{disarray}$$
which is bounded by a constant independent of $t$ by Remark
\ref{lm2}, and as $\rho <v$, we have $\lim_{t\rightarrow\infty} e^{
(\rho-v)t}=0$. Therefore:
\begin{eqnarray}
\lim_{t\rightarrow\infty}\mathbb{E}^{\updownarrow} (c_{t})=0\
.\label{ct}
\end{eqnarray}

 $\bullet$ Now we consider $d_{t}$. As shown in \cite{berrou}
 with $B (r,j)=\{i\geq 2:
 \Pi_{i} (t)\subseteq\pi_{j} (r)\cap\Pi_{1} (r-)\}$, we get, for
 every $r\in[0,t]$ and $j\geq 2$, conditionally on
 $r\in\mathcal{D}_{1},\ \  \Pi_{1} (r-)$ and $\pi_{j} (r)$, the
 partition $ (\Pi_{i} (t): i \in B (r,j))$ can be written in the form
$\tilde\Pi^{(j)} (t-r)_{|\pi_{j} (r)\cap\Pi_{1} (r-)}$. Here $(\tilde\Pi^{(j)})_{j\in\mathbb{N}}$
is a family of i.i.d. homogeneous fragmentations distributed as $\Pi$ under
$\mathbb{P}$ and independent of the sigma-field $\mathcal{G}_{1}
(t)$. As a consequence:
$$\underset{i\geq 2}{\cup}\Pi_{i} (t)=\underset{j\geq2}{\cup}\underset{r\in[0,t]\cap
\mathcal{D}_{1}}{\cup}\tilde\Pi^{(j)} (t-r)_{|\pi_{j} (r)\cap\Pi_{1}
(r-)}.$$Moreover $|\pi_{j} (r)||\Pi_{1} (r-)|$ is $\mathcal{G}_{1}
(t)$
 measurable, and we have that for all $i\in\mathbb N$ \[|\tilde{\Pi_{i}}^{(j)} (t-r)_{|\pi_{j} (r)\cap\Pi_{1} (r-)}|
 =|\tilde{\Pi_{i}}^{(j)} (t-r)||\pi_{j} (r)||\Pi_{1} (r-)|\] so that we get:
$$\begin{disarray}{rcl}
\mathbb{E}^{\updownarrow} (d_{t}|\mathcal{G}_{1} (t))&\leq &
 \frac{e^{\rho t}}{h (0)}C_{8}
\sum_{r\in[0,t]\cap \mathcal{D}_{1} }\sum_{j=2}^{\infty}|\pi_{j}
(r)||\Pi_{1}^{\dag} (r-)|\nbOne_{\{a\leq |\pi_{j}
(r)||\Pi_{1}^{\dag} (r-)|e^{vr}\leq b\}}\\&&
\sum_{i=1}^{\infty}\mathbb{E}^{\updownarrow} \left(
 |\tilde{\Pi_{i}}^{(j)} (t-r)\nbOne_{\{a\leq |\tilde{\Pi_{i}}^{(j)} (t^{'}-r)|e^{v(t^{'}-r)}|\pi_{j}
(r)||\Pi_{1}^{\dag} (r-)|e^{vr}\leq b \ \forall t^{'}\in[r,t]\}}\
|\mathcal{G}_{1} (t)\right)\\&\leq &
 \frac{e^{\rho t}}{h (0)}C_{8}
\sum_{r\in[0,t]\cap \mathcal{D}_{1} }\sum_{j=2}^{\infty}|\pi_{j}
(r)||\Pi_{1}^{\dag} (r-)|\nbOne_{\{a\leq |\pi_{j}
(r)||\Pi_{1}^{\dag} (r-)|e^{vr}\leq b\}}\\&&
\sum_{i=1}^{\infty}\mathbb{E}^{\updownarrow} \left(
 |\tilde{\Pi_{i}}^{(j)} (t-r)\nbOne_{\{a/b\leq |\tilde{\Pi_{i}}^{(j)}(t^{'}-r)|e^{v(t^{'}-r)}\leq b/a   \ \forall
t^{'}\in[r,t]\}}\ |\mathcal{G}_{1} (t)\right),
\end{disarray}$$
 with $C_{8}$ the maximum of
$h (t)$ on the compact $[\log (a),\log (b)]$.  As $\tilde\Pi$ is
independent of the sigma-field $\mathcal{G}_{1} (t)$, $\tilde\Pi$
has the same distribution under $\mathbb{P}$ as under
$\mathbb{P}^{\updownarrow}$. $\tilde\Pi^{(j)}$  is also distributed as
$\Pi$ under $\mathbb{P}$. Thus,
$$\begin{disarray}{l}\left. \sum_{i=1}^{\infty}\mathbb{E}^{\updownarrow} \left (
 |\tilde{\Pi_{i}}^{(j)} (t-r)\nbOne_{\{a/b\leq |\tilde{\Pi_{i}}^{(j)} (t^{'}-r)|e^{v(t^{'}-r)}\leq b/a   \ \forall
t^{'}\in[r,t]\}}|\ \right|\mathcal{G}_{1}
(t)\right)\\=\sum_{i=1}^{\infty}\mathbb{E}^{\updownarrow}  \left (
 |\tilde{\Pi_{i}}^{(j)} (t-r)|\nbOne_{\{a/b\leq |\tilde{\Pi_{i}}^{(j)} (t^{'}-r)|e^{v(t^{'}-r)}\leq b/a   \ \forall
t^{'}\in[r,t]\}}\ \right).\end{disarray}$$

 Now we have by
size-biased sampling:
$$\begin{disarray}{l}
\sum_{i=1}^{\infty} \mathbb{E}\left (e^{\rho (t-r)}|\Pi_{i}
(t-r)|\nbOne_{\{a/b\leq |\Pi_{i}
(t^{'}-r)|e^{v(t^{'}-r)}\leq b/a   \ \forall
t^{'}\in[r,t]\}}\right)\\= \mathbb{E}\left (\ e^{\rho (t-r)}
\mathbf{1}_{\{t-r<\inf\{s:\ |\Pi_{1} (s)| \notin \ (\frac{a}{b}
e^{-vs},\frac{b}{a}e^{-vs})\}\}}\right)\\ \\
= \mathbf{E}_{\log (1/a)}\left (\ e^{\rho (t-r)}
\mathbf{1}_{\{T_{2\log (b/a)}>t-r\}\}}\right) ,
\end{disarray}$$
as $\rho_{.}$ is decreasing $\rho_{2\log (b/a)}\leq \rho$, thus
$$\begin{disarray}{l}
\sum_{i=1}^{\infty} \left.e^{\rho (t-r)}\mathbb{E}^{\updownarrow} \left (
 |\tilde{\Pi_{i}}^{(j)} (t-r)\nbOne_{\{a/b\leq |\tilde{\Pi_{i}} (t^{'}-r)|e^{v(t^{'}-r)}\leq b/a   \ \forall
t^{'}\in[r,t]\}}|\ \right|\mathcal{G}_{1}
(t)\right)\\\leq \mathbf{E}_{\log (1/a)}\left (\
e^{\rho_{2\log (b/a)} (t-r)} \mathbf{1}_{\{T_{2\log
(b/a)}>t-r\}\}}\right) .
\end{disarray}$$
Therefore  with $F_{x}^{'} (t):=\mathbf{E}_{x}\left (\
e^{\rho_{2\log (b/a)} t} \mathbf{1}_{\{T_{2\log
(b/a)}>t\}\}}\right)$  and  since $|\Pi_{1}^{\dag}
(r-)|=ae^{Y_{r-}-vr}\nbOne_{\{r-<T\}}$ under $\mathbb{P}_{\log
(1/a)}^{\updownarrow}$, we get:

$$\begin{disarray}{rcl}
\mathbb{E}^{\updownarrow} (d_{t}|\mathcal{G}_{1} (t))&\leq& a
\frac{e^{ (\rho-v) r}C_{8}}{h (0)} \sum_{r\in[0,t]\cap
\mathcal{D}_{1} }\sum_{j=2}^{\infty}|\pi_{j}
(r)|e^{Y_{r-}}\sup_{x\in[0, \log (b/a)]}\sup_{t\geq 0}|F_{x}^{'}
(t)|.
\end{disarray}$$

Moreover  we have by definition
$e^{Y_{r-}}\nbOne_{\{r-<T\}}\leq b/a$. We let
\[C_{9}:=\sup_{x\in[0, \log (b/a)]}\sup_{t\geq 0}|F_{x}^{'} (t)|bC_{8}/ h
(0)\]according to Remark~\ref{lm2} we have $C_{9}<\infty$. Thus
$$\begin{disarray}{rcl}
\mathbb{E}^{\updownarrow} (d_{t}|\mathcal{G}_{1} (t))&\leq
&C_{9}\sum_{r\in[0,t]\cap \mathcal{D}_{1}  }\sum_{j=2}^{\infty}e^{
(\rho-v) r}|\pi_{j} (r)|.
\end{disarray}$$

Under $\mathbb{P}$, the $\mathcal{G}_{0} (t-)$ predictable
compensator of \[A_{t}:=\sum_{r\in[0,t]\cap \mathcal{D}_{1} }
\sum_{j=2}^{\infty}e^{ (\rho-v) r}|\pi_{j} (r)|\] is
\[N_{t}:=\int_{0}^{t}dr \int_{\mathcal{P}}\mu (ds) e^{ (\rho-v)r}
\sum_{j=2}^{\infty}|\pi_{j}| .\] Additionally
\[\int_{\mathcal{P}}\mu (ds)
\sum_{j=2}^{\infty}|\pi_{j}| = \int_{\mathcal{S}^{*}}\nu (ds)
 \sum_{i=1}^{\infty}s_{i}\left[(\sum_{j=1}^{\infty}s_{j})-s_{i}\right].\]

As   $\sum_{j=1}^{\infty} s_{j}=1$ $\nu - a.s.$, we achieve:

\[\int_{\mathcal{P}}\mu (ds)
\sum_{j=2}^{\infty}|\pi_{j}| \leq \int_{\mathcal{S}^{*}}\nu (ds)
 2(1-s_{1}),\]
which is finite by (\ref{mesuredelevy}). Moreover as $\rho <v$, the
term $e^{(\rho-v)r}$ is integrable on $[0,\infty)$, so that we have
$\lim_{t\rightarrow\infty}N_{t}<\infty$.

 As both $X_{t}:=A_{t}-N_{t}$  and
$M_{t}$ are martingales,  by Theorem 4.50 of \cite{jacshi}, we get
that  \[XM-[X,M]\ \ \ \  \text{is a local martingale}.\]

A sequence $(\tau_{n}=(T(m,n))_{m\in\mathbb{N}})_{n\in\mathbb{N}}$
of adapted subdivisions is called a Riemann sequence if
$\sup_{m\in\mathbb{N}}[T(m+1,n)\wedge t-T(m,n)\wedge t]\rightarrow
0$ for all $t\in\mathbb{R}_{+}$.  By Theorem 4.47 of \cite{jacshi},
for any Riemann sequence
$\{\tau_{n}=(T(m,n))_{m\in\mathbb{N}}\}_{n\in\mathbb{N}}$ of adapted
subdivisions, the processes $S_{\tau_{n}}(X,M)$ defined by
\[S_{\tau_{n}}(X,M)_{t}:=\sum_{m\in\mathbb{N}}(X_{T(m+1,n)\wedge t}-X_{T(m,n)\wedge t})(M_{T(m+1,n)\wedge t}-M_{T(m,n)\wedge t})\]
converge to the process $[X,M]$, in measure, uniformly on every
compact interval.

We will now bound $S_{\tau_{n}}(X,M)_{t}$ uniformly in $t$. As
\begin{equation}S_{\tau_{n}}(X,M)_{t}\leq
\sup_{l\in\mathbb{N}}|M_{T(l+1,n)\wedge t}-M_{T(l,n)\wedge
t}|\sum_{m\in\mathbb{N}}|X_{T(m+1,n)\wedge t}-X_{T(m,n)\wedge
t}|\label{eqtau}\end{equation} we will first focus on
$\sum_{m\in\mathbb{N}}|X_{T(m+1,n)\wedge t}-X_{T(m,n)\wedge t}|$:

$$\begin{disarray}{l}\sum_{m\in\mathbb{N}}|X_{(m+1)/n\wedge
t}-X_{T(m,n)\wedge t}|\\\leq \sum_{m\in\mathbb{N}}\left(\sum_{r\in[T(m,n)\wedge
t,T(m+1,n)\wedge t]\cap \mathcal{D}_{1}  }\sum_{j=2}^{\infty}e^{
(\rho-v) r}|\pi_{j} (r)| +\int_{T(m,n)\wedge t}^{T(m+1,n)\wedge t}dr
\int_{\mathcal{P}}\mu (ds) e^{ (\rho-v)r}
\sum_{j=2}^{\infty}|\pi_{j}|\right)\\
\leq \sum_{r\in[0,t]\cap \mathcal{D}_{1}  }\sum_{j=2}^{\infty}e^{
(\rho-v) r}|\pi_{j} (r)|+\int_{0}^{\infty}dr \int_{\mathcal{P}}\mu
(ds) e^{ (\rho-v)r} \sum_{j=2}^{\infty}|\pi_{j}|.\end{disarray}$$

Therefore by the previous study of $A_{t}$ and $N_{t}$ we get that
there exist $C_{10}<\infty$ independent of $t$ such that:
\[\lim_{n\rightarrow\infty}\mathbb{E}\left( \sum_{m\in\mathbb{N}}|X_{T(m+1,n)\wedge t}-X_{T(m,n)\wedge
t}|\right)\leq C_{10}\ \ \  \text{for all } \ t.\] Moreover
\[\lim_{n\rightarrow\infty}\sup_{l\in\mathbb{N}}|M_{T(l+1,n)\wedge t}-M_{T(l,n)\wedge
t}|\leq \sup_{r\leq t}|M_{r}-M_{r-}|\] is a.s. bounded by $C^{'}$
(see remark \ref{rkmt}) independently of $t$. Consequently by
(\ref{eqtau})

\[\lim_{t\rightarrow\infty}\mathbb{E}([X,M]_{t})<\infty .\]
Thus as $XM-[X,M]$ is a local martingale, we get  that
$\lim_{t\rightarrow\infty}\mathbb{E}^{\updownarrow} (d_{t})<\infty\
.$

Finally according to ~ (\ref{ct}), we get
$$\lim_{t\rightarrow\infty}\mathbb{E}^{\updownarrow} (M_{t})
=\lim_{t\rightarrow\infty}\mathbb{E}^{\updownarrow}
(\mathbb{E}^{\updownarrow} (d_{t}+c_{t}|\mathcal{G}_{1} (t)))<\infty
.
$$

\bigskip

\noindent{\bf Acknowledgments: } I wish to thank J. Bertoin for
his help, and especially for his patience.

I am very grateful to A. Kyprianou, who kindly pointed out an
original mistake in the proof of Theorem \ref{theo2}.1. and for his
valuable comments and advices.

I also wish to thank an anonymous referee of an earlier draft for  various interesting suggestions.

\end{document}